\documentclass{article}
\usepackage{graphicx} 
\usepackage[utf8]{inputenc}
\usepackage{array} 
\usepackage{amsfonts}
\usepackage{amsmath}
\usepackage{amsthm}

\usepackage{amssymb}
\usepackage{xfrac}
\usepackage[margin=1in]{geometry}
\usepackage{subcaption} 

\usepackage{biblatex}
\usepackage{hyperref}

\theoremstyle{definition}

\newcommand{\R}{\mathbb{R}}

\newcommand{\domain}{\mathcal{D}}
\DeclareMathOperator{\vv}{\textbf{v}}
\DeclareMathOperator{\vu}{\textbf{u}}

\DeclareMathOperator{\vn}{\textbf{n}}
\DeclareMathOperator{\f}{\textbf{f}}
\DeclareMathOperator{\g}{\textbf{g}}
\DeclareMathOperator{\F}{\mathcal{F}}
\newcommand{\T}{\mathbf{T}}
\newcommand{\ess}{\mathbf{S}}

\newcommand{\A}{\mathbf{A}}

\newcommand{\B}{\mathbf{B}}

\newcommand{\W}{\mathbf{W}}
\newcommand{\D}{\mathbf{D}}

\newcommand{\Sig}{\boldsymbol{\Sigma}}

\NewDocumentCommand{\del}{e{_^}}{%
  \mathop{}\!
  \nabla
  \IfValueT{#1}{_{\!#1}}
  \IfValueT{#2}{^{#2}}
}

\newcommand{\G}{\mathcal{G}}

\DeclareMathOperator{\Wi}{Wi}
\DeclareMathOperator{\Rey}{Re}

\title{A Finite Element Implementation of the SRTD Algorithm for an Oldroyd 3-Parameter Viscoelastic Fluid Model}
\author{Christian Austin \\ Sara Pollock \\ L. Ridgway Scott}
\date{\today}

\begin{document}

\maketitle

\begin{abstract}
In this paper, we discuss a finite element implementation of the SRTD algorithm described by Girault and Scott for the steady-state case of a certain 3-parameter subset of the Oldroyd models. We compare it to the well-known EVSS method, which, though originally described for the upper-convected Maxwell model, can easily accommodate the Oldroyd 3-parameter model. We obtain numerical results for both methods on two benchmark problems: the lid-driven cavity problem and the journal-bearing, or eccentric rotating cylinders, problem. We find that the resulting finite element implementation of SRTD is stable with respect to mesh refinement and is generally faster than EVSS, though is not capable of reaching as high a Weissenberg number as EVSS.
\end{abstract}

\section{Introduction}\label{sec1:intro}
Non-Newtonian fluids abound in both the natural sciences and engineering and include a wide array of materials such as lubricants, oils, polymer melts, shampoo, tar, blood, quicksand, and wet cement \cite{rheology-macosko}. Even the flow of food emulsions such as mayonnaise and the gradual creeping flow of ice can be predicted by appropriate non-Newtonian models \cite{Dunn_Rajagopal_g2}. The modeling and simulation of such fluids therefore has wide applicability. 

Viscoelastic fluids are a special subset of non-Newtonian fluids that are capable of storing and releasing elastic energy. We consider the steady-state case of a certain 3-parameter model for viscoelastic fluids which is a subset of the Oldroyd models from \cite{oldroyd1958}. The 3 parameters are the kinematic viscosity $\eta_{0}$ (we will assume a constant density of $\rho=1$), the relaxation time $\lambda_{1}$, and another rheological parameter $\mu_{1}$ associated with the slip parameter. This parameter subset includes the upper-convected, lower-convected, and corotational Maxwell models. 

In certain other viscoelastic fluid models, there exists another rheological parameter called the retardation time $\lambda_{2}$. When this value is is strictly positive, an explicit dissipation term or Laplacian appears in the momentum equation, making it explicitly elliptic and lending itself well to finite element approximation. Notably, the Maxwell models we study here have zero retardation time, $\lambda_{2}=0$, meaning the momentum equation lacks explicit dissipation. Various methods have been introduced in the literature which manipulate the resulting system of equations for a Maxwell-like fluid to introduce a Laplacian term into the momentum equation; some notable historical examples include the Explicitly Elliptic Momentum Equation (EEME) method \cite{eeme} and the Elastic-Viscous Split Stress (EVSS) method \cite{evss}. Recently, a new formulation, termed Selective Replacement of Tensor Divergence, or SRTD, has been developed by Girault and Scott in \cite{scott_wellposed} for introducing an explicit dissipation term into the momentum equation. Their work is based on the earlier work of Renardy \cite{renardy_iter} for the upper-convected Maxwell model and, like Renardy's method, includes an iterative scheme for solving the system. The SRTD iteration was used in \cite{scott_summary, scott_wellposed} to establish well-posedness of Maxwell-like models in Sobolev spaces, but this paper is the first to examine it computationally. By contrast, EVSS has received significant computational attention, but so far does not have a rigorous basis in terms of convergence in any spaces. Here, we implement the SRTD iterative scheme in a finite element context and discuss the results of numerical experiments.

\subsection{Notation}\label{subsec1:1_notation}
For consistent notation, we declare now that the dot product $\cdot$ will refer to either usual matrix multiplication or matrix-vector multiplication, with the choice of the two being clear from the context. That is, for $2$-tensors (matrices) $\A$ and $\B$ and for a $1$-tensor (vector) $\vu$, using Einstein notation, we have
\[(\A\cdot\B)_{ij} = \A_{ik}\B_{kj}, \quad\text{ and }\quad (\B\cdot\vu)_{i} = \B_{ij}\vu_{j}.\]
We let $\A:\B$ denote tensor contraction, 
\[\A:\B = \A_{ij}\B_{ij}.\]
Letting an index after a comma denote differentiation with respect to that spatial variable, we say the gradient of a $1$-tensor $\vu$ has entries
\[\quad(\del\vu)_{ij}= \frac{\partial \vu_{i}}{\partial x_{j}} = \vu_{i,j} \]
and that the divergence of a $2$-tensor $\B$ has entries
\[(\del\cdot\B)_{i}= \B_{ij,j}.\]
This choice of the indexing for the gradient and divergence means that the Laplacian $\Delta\vu$ of a vector field is given by the divergence of the gradient,
\[(\Delta\vu)_{i} = \vu_{i,jj} = (\del\cdot\del\vu)_{i}\]
but that the advective derivative $\vu\cdot\del$ of a tensor quantity $\B$ is ``backwards," 
\[(\vu\cdot\del\B)_{ij} := \B_{ij,k}\vu_{k} = (\del\B\cdot\vu)_{ij}. \]
For clarity, we will declare now that $\vu\cdot\del\B$ always refers to the advective derivative of $\B$. 

For a domain $\domain\subset\R^{d}$ with spatial dimension $d\in\{2,3\}$, we will use $W_{q}^{s}(\domain)$ to indicate the the standard Sobolev space, or the space of all functions on $\domain$ whose weak derivatives up to order $s$ all live in $L_{q}(\domain)$. We use $W_{q}^{s}(\domain)^{d}$ and $W_{q}^{s}(\domain)^{d^{2}}$ for vector and $2$-tensor--valued functions, respectively, whose components all live in $W_{q}^{s}(\domain)$. 

We use $\domain_{h}$ to denote a simplicial (triangular or tetrahedral) discretization of the domain $\domain$ with characteristic mesh size $h$. Given a discretization $\domain_{h}$ of our domain, we let $P^{k}(\domain_{h})$ denote continuous piecewise polynomials of degree at most $k$ on $\domain_{h}$. Likewise, let $P^{k}(\domain_{h})^{d}$ and $P^{k}(\domain_{h})^{d^{2}}$ denote vector and tensor-valued functions, respectively, whose components all live in $P^{k}(\domain_{h})$. A subscript $0$, for example, $P^{k}_{0}(\domain_{h})^{d}$, refers to functions in the respective spaces which have boundary trace zero. 

\subsection{Equations of Motion}\label{subsec1:2_equations_of_motion}
All steady or time-invariant fluid motion can be described via the momentum equation, which is an expression of Newton's second law written in a way convenient for fluids, 
\begin{equation}
    \vu\cdot\del\vu = -\del p +\del\cdot\T + \f,
    \label{eq:momentum}
\end{equation}
where $\vu$ is the fluid velocity, $p$ is the pressure, $\T$ is the extra stress tensor, and $\f$ represents any external body forces. The momentum equation \eqref{eq:momentum} is accompanied by the continuity equation, which is an expression of the law of conservation of mass. For an incompressible fluid of constant density, the continuity equation takes the form 
\begin{equation}
    \del\cdot\vu = 0.
    \label{eq:continuity}
\end{equation}
To complete the model, one must describe how the stress tensor $\T$ relates to the fluid velocity via a constitutive equation. For a Newtonian fluid, the constitutive equation reads
\begin{equation}
    \T = 2\eta_{0}\D,
    \label{eq:newtonian}
\end{equation}
where  $\eta_{0}$ is the kinematic viscosity, and $\D$ is the fluid rate of deformation or rate of strain tensor,
\begin{equation}
    \D = \D(\vu) = \frac{1}{2}(\del\vu+\del\vu^{t}).
    \label{eq:deformation_tensor}
\end{equation}
Taken with the incompressibility assumption \eqref{eq:continuity}, this implies $\del\cdot\T = \eta_{0}\Delta\vu$. This can be substituted directly into \eqref{eq:momentum}, yielding the familiar steady-state Navier-Stokes equations, 
\begin{equation}
    -\eta\Delta\vu + \vu\cdot\del\vu + \del p= \f,
    \label{eq:nse}
\end{equation}
which must be solved together with \eqref{eq:continuity}.
For viscoelastic and other non-Newtonian fluids, the constitutive equation used in place of \eqref{eq:newtonian} is more complicated. 

\subsection{Objective Rates}\label{subsec1:3_objective_rates}
The constitutive equations defining many viscoelastic fluid models are often expressed using convected derivatives or objective rates. These are material time derivatives with additional terms added to describe a quantity which stretches or compresses with the fluid and to make the entire quantity frame invariant \cite{cara_guillen_ortega_modeling_oldroyd, reddy-mechanics}. The upper-convected time derivative (UCTD) of a $2$-tensor $\B$, denoted $\overset{\bigtriangledown}{\B}$, is such an objective rate, and is given by
\begin{equation}
    \overset{\bigtriangledown}{\B} := \frac{\partial\B}{\partial t} + \vu\cdot\del\B - (\del \vu)\cdot \B - \B\cdot(\del \vu)^{t},
    \label{eq:uctd}
\end{equation}
where $\vu$ is the fluid velocity. Another is the lower-convected time derivative (LCTD), denoted $\overset{\triangle}{\B}$, 
\begin{equation}
    \overset{\triangle}{\B} := \frac{\partial\B}{\partial t} + \vu\cdot\del\B + (\del \vu)^{t}\cdot \B + \B\cdot(\del \vu).
    \label{eq:lctd}
\end{equation}
Any convex combination of the UCTD and LCTD results in a valid objective rate. This range of objective rates is often expressed via a dimensionless parameter $a\in[-1,1]$, sometimes called the \textit{slip parameter} \cite{rheology-morrison}. We let $\G(\cdot,a)$ denote such a general convected derivative, 
\begin{equation}
    \G(\B, a) := \frac{\partial\B}{\partial t} + \vu\cdot\del\B + \W\cdot\B - \B\cdot\W - a (\D\cdot\B + \B\cdot\D),
    \label{eq:general_ctd}
\end{equation}
where $\D$ is the fluid rate of deformation tensor defined in \eqref{eq:deformation_tensor}, and $\W$ is the fluid spin tensor, 
\begin{equation}
    \W = \W(\vu) = \frac{1}{2}(\del\vu^{t} - \del\vu).
    \label{eq:spin_tensor}
\end{equation}
Note that $\G(\B, -1)$ is the LCTD of $\B$, and $\G(\B, 1)$ is the UCTD. Another special case is when $a=0$, and $\G(\B, 0)$ is often called the \textit{corotational derivative}, 
\begin{equation}
    \G(\B, 0) = \frac{\partial\B}{\partial t} + \vu\cdot\del\B + \W\cdot\B - \B\cdot\W.
    \label{eq:corot_ctd}
\end{equation}

\subsection{Viscoelastic Fluid Models}\label{subsec1:4_viscoelastic_models}

 One of the simplest viscoelastic fluid models is the upper-convected Maxwell (UCM) model, first described by Oldroyd \cite{oldroyd1950} in 1950 and based on a model of Maxwell for describing the dynamics of gases \cite{maxwell}. The constitutive equation of the UCM model reads
\begin{equation}
    \T + \lambda_{1}\overset{\bigtriangledown}{\T} = 2\eta_{0}\D,
    \label{eq:ucm}
\end{equation}
where $\lambda_{1}$ is the fluid relaxation time, $\eta_{0}$ the kinematic viscosity, and $\overset{\bigtriangledown}{\T}$ is the UCTD of $\T$, given in \eqref{eq:uctd}. An immediate generalization can be obtained by replacing the UCTD of $\T$ in the UCM model \eqref{eq:ucm} with an arbitrary convected derivative $\G(\T, a)$, 
\begin{equation}
    \T + \lambda_{1}\G(\T, a) = 2\eta_{0}\D,
    \label{eq:general_maxwell}
\end{equation}
which, when fully expanded via \eqref{eq:general_ctd}, in the steady-state case takes the form
\begin{equation}
    \T + \lambda_{1}\bigl( \vu\cdot\del\T + \W\cdot\T - \T\cdot\W - a (\D\cdot\T + \T\cdot\D) \bigr) = 2\eta_{0}\D.
    \label{eq:general_maxwell_expanded}
\end{equation}
We refer to the viscoelastic model whose constitutive equation is described by \eqref{eq:general_maxwell} as a generalized Maxwell model. This model can also be viewed as a 3-parameter subset of an 8-parameter model described by Oldroyd \cite{oldroyd1958} in 1958. Under the Oldroyd 8-parameter scheme, the relaxation time $\lambda_{1}$ and kinematic viscosity $\eta_{0}$ play their usual roles, while the dimensionless slip parameter $a$ is replaced with the parameter $\mu_{1}=a\lambda_{1}$, which has dimension of time. The remaining 5 parameters of the model in \cite{oldroyd1958} are set to zero, and the constitutive equation takes the equivalent form 
\begin{equation}
    \T + \lambda_{1}( \vu\cdot\del\T + \W\cdot\T - \T\cdot\W) - \mu_{1} (\D\cdot\T + \T\cdot\D) = 2\eta_{0}\D.
    \label{eq:oldroyd3}
\end{equation}
Adding and subtracting $\lambda_{1}(\D\cdot\T + \T\cdot\D)$, it can also be rewritten as 
\begin{equation}
    \T + \lambda_{1}\bigl(\vu\cdot\del\T - (\del\vu)\cdot\T - \T\cdot(\del\vu)^{t}\bigr) + (\lambda_{1}-\mu_{1})(\D\cdot\T + \T\cdot\D) = 2\eta_{0}\D,
    \label{eq:oldroyd3_alt}
\end{equation}
which is how it is commonly written in \cite{scott_wellposed}. In this context, we refer to the model as the \textit{Oldroyd 3 parameter subset model} (O3). This is the primary lens through which we will view this model, so we will generally use $\mu_{1}$ rather than $a$ to characterize the convected derivative. Note that when $\mu_{1}=\lambda_{1}$, we indeed recover the UCM model. 

The immediate difficulty surrounding many viscoelastic models, including O3, when compared to the Newtonian or Navier-Stokes case \eqref{eq:nse} is that the divergence of the stress tensor $\del\cdot\T$ cannot be easily isolated from the divergence of the constitutive equation \eqref{eq:general_maxwell} and substituted directly into \eqref{eq:momentum}.  Therefore, a coupled system of equations involving \eqref{eq:momentum}, \eqref{eq:continuity}, and \eqref{eq:general_maxwell}, \eqref{eq:general_maxwell_expanded}, \eqref{eq:oldroyd3}, or \eqref{eq:oldroyd3_alt} must be solved. The full system of equations for the steady flow of an incompressible fluid governed by the O3 model reads 
\begin{equation}
    \begin{aligned}
        \vu\cdot\del\vu + \del p &= \del\cdot\T + \f, \\
        \del\cdot\vu &= 0, \\
        \T + \lambda_{1}( \vu\cdot\del\T + \W\cdot\T - \T\cdot\W) - \mu_{1} (\D\cdot\T + \T\cdot\D) &= 2\eta_{0}\D.
    \end{aligned}
    \label{eq:oldroyd3_system}
\end{equation}

Additional difficulties with solving \eqref{eq:oldroyd3_system} which are not present in the Newtonian or Navier-Stokes case come from the fact that the momentum equation does not contain an explicit Laplacian term. In models like the Johnson-Segalman model \cite{johnson-segalman}, and its most well-known special case, the Oldroyd B model \cite{oldroyd1950}, the extra stress tensor $\T$ can be decomposed into a solvent part and a polymer part, $\T = \T_{s}+\T_{p}$, where $\T_{s}=2\eta_{s}\D$ represents a contribution from a Newtonian solvent, and $\T_{p}$ represents a contribution from a viscoelastic polymer modeled by the O3 constitutive equation \eqref{eq:general_maxwell}. In such models where there is a Newtonian solvent contribution, the divergence of the stress tensor introduces a Laplacian operator into the momentum equation with coefficient equal to the solvent viscosity, as
\[\del\cdot\T = \del\cdot\T_{s} + \del\cdot\T_{p} = \eta_{s}\Delta\vu + \del\cdot\T_{p}.\]
Such a term makes the momentum equation explicitly elliptic and aids in stability of numerical schemes. The O3 model, however, is precisely the special case of the Johnson-Segalman model that has no Newtonian solvent contribution.  Instabilities can arise when there is a lack of a Newtonian solvent contribution, as is the case for the UCM and O3 models \cite{alves_review, scott_wellposed, change_of_type_joseph_renardy_saut}.

The remainder of the paper is organized as follows. In Section \ref{sec2:methods}, we discuss two methods of introducing a Laplacian term into the momentum equation and thus helping to stabilize the system \eqref{eq:oldroyd3_system}, making it suitable for a finite element approximation. The first is the Elastic-Viscous Split Stress (EVSS) method \cite{evss}, which was first introduced in 1990 by Rajagopalan, Armstrong, and Brown and will serve as our control or baseline method. The second is the newer SRTD formulation \cite{scott_wellposed}, introduced by Girault and Scott in 2017. In Section \ref{sec3:test_problems}, we introduce the physical flow problems, the journal-bearing or eccentric-rotating cylinders problem, and the lid-driven cavity problem, on which we will perform our numerical experiments. Section \ref{sec4:numerical_results} includes the results of numerical experiments on the convergence and stability of the SRTD method, and we include some comparisons with EVSS. 

\section{Formulations and Methods}\label{sec2:methods}

\subsection{EVSS}\label{subsec2:1_EVSS}
The elastic-viscous split stress (EVSS) formulation was first described in 1990 by Rajagopalan, Armstrong, and Brown \cite{evss}. In their paper, the EVSS formulation was only given explicitly for the UCM and Giesekus models, although the authors remark that the EVSS formulation can be ``easily extended to a wide variety of constitutive equations, with or without a solvent viscosity." This is indeed the case, and we do so here for the O3 model. 

The idea behind EVSS in general is to introduce a new variable, called the \textit{elastic stress tensor} $\Sig$, to isolate the viscous contribution from the constitutive equation. This explicitly introduces a Laplacian term into the momentum equation similar to how the solvent contribution does for Johnson-Segalman fluids. For the O3 model, we define the elastic stress tensor $\Sig$ to be
\begin{equation}
    \Sig := \T - 2\eta_{0}\D.
    \label{eq:elastic_stress}
\end{equation}
A constitutive equation for $\Sig$ can be easily derived from \eqref{eq:general_maxwell}, 
\begin{equation}
    \Sig + \lambda_{1}\G(\Sig + 2\eta_{0}\D,a) = 0, 
    \label{eq:evss_general_maxwell}
\end{equation}
and, rewriting the constitutive equation for $\Sig$ like \eqref{eq:oldroyd3_alt}, the full system of equations for the EVSS formulation of the O3 model is given by 
\begin{equation}
    \begin{aligned}
        -\eta_{0}\Delta\vu + \vu\cdot\del\vu + \del p &= \del\cdot\Sig + \f, \\
        \del\cdot\vu &= 0, \\
        \Sig + \lambda_{1}\Bigl(\vu\cdot\del(\Sig + 2\eta_{0}\D) - (\del\vu)\cdot(\Sig + 2\eta_{0}\D) - (\Sig + 2\eta_{0}\D)\cdot(\del\vu)^{t}\Bigr) &\\
        + (\lambda_{1}-\mu_{1})\bigl(\D\cdot(\Sig + 2\eta_{0}\D) + (\Sig + 2\eta_{0}\D)\cdot\D\bigr) &= 0.
    \end{aligned}
    \label{eq:evss_formulation}
\end{equation}
Notice the inclusion of a Laplacian term in the momentum equation.

\subsubsection{EVSS Finite Element Implementation}\label{subsubsec2:1_1_EVSS_FE}

The immediate drawback of applying EVSS to the O3 model, and in many other cases where EVSS can be applied, is the introduction of a higher derivative of $\vu$ when the convected derivative of $\D$ is taken in \eqref{eq:evss_general_maxwell}, making \eqref{eq:evss_formulation} somewhat unsuitable for direct finite element approximation. Rajagopalan et al. remedy this by taking the rate of strain tensor $\D$ to be a separate variable and enforcing its definition $\D = \D(\vu)$ via an $L_{2}$ projection, manifesting as an additional equation in the finite element formulation. Since $\D$ must be symmetric and traceless, they remark that this introduces only $2$ new equations into the overall mixed formulation for the 2-dimensional case, or $5$ for the 3-dimensional case. The EVSS formulation of the O3 model, suitable for finite element approximation, is given by 

\begin{equation}
    \begin{aligned}
        -\eta_{0}\Delta\vu + \vu\cdot\del\vu + \del p &= \del\cdot\Sig + \f, \\
        \del\cdot\vu &= 0, \\
        \Sig + \lambda_{1}\Bigl(\vu\cdot\del(\Sig + 2\eta_{0}\D) - (\del\vu)\cdot(\Sig + 2\eta_{0}\D) - (\Sig + 2\eta_{0}\D)\cdot(\del\vu)^{t}\Bigr) &\\
        +(\lambda_{1}-\mu_{1}) \Bigl(\D(\vu) \cdot(\Sig + 2\eta_{0}\D) + (\Sig + 2\eta_{0}\D)\cdot\D(\vu)\Bigr) &= 0, \\ 
        \D &= \frac{1}{2}\Bigl(\del\vu + (\del\vu)^{t}\Bigr),
    \end{aligned}
    \label{eq:evss_system_D}
\end{equation}
and the stress tensor $\T$ can be optionally recovered from \eqref{eq:elastic_stress} with $\T = \Sig + 2\eta_{0}\D$. Notice that $\D$ is only treated as a separate variable when it is subject to a convected derivative (in order to avoid taking higher derivatives of $\vu$ explicitly). When this is not the case, such as when $\D$ appears as a coefficient in the definition of a general convected derivative, $\D(\vu) = \frac{1}{2}(\del\vu + \del\vu^{t})$ is used. 

Since the constitutive equation \eqref{eq:evss_general_maxwell} is an advection equation for for $\Sig$, Rajagopalan et al. recommend using streamline-upwind Petrov-Galerkin (SUPG) stabilization \cite{supg_brooks_hughes}. SUPG is often used to stabilize the finite element solution to advection dominated flows, and preexisting finite element codes can often be easily adjusted to include SUPG. 

As we are solving the fully mixed formulation \eqref{eq:evss_system_D} using FEM, we would like to ensure our mixed space satisfies the necessary compatibility conditions for wellposedness. For viscoelastic models like the UCM, however, the situation is complicated, and most of what is known is based on the 3-field $(\vu, p, \T)$ formulation of the Stokes equations or other linearizations of the equations \eqref{eq:oldroyd3_system}; the situation is summarized nicely by Owens and Philips in Section 7.3 of \cite{computational_rheology}. Baranger and Sandri in \cite{Baranger_sandri_3_field_stokes}, for instance, study the 3-field $(\vu, p, \T)$ Stokesian limit of the UCM model and find that, when there is a Newtonian contribution to the momentum equation (that is, a Laplacian term is still present in the momentum equation), no additional compatibility conditions are necessary other than the usual LBB conditions for the velocity-pressure pair. Owens and Philips \cite{computational_rheology} further note that when an EVSS-like change of variables is used to produce a Laplacian in the momentum equation, no further comppatibility conditions are needed than the usual LBB conditions in the Stokesian limit. 

Rajagopalan et al. in their paper describing EVSS \cite{evss} use a quadrilateral mesh to discretize their domain and recommend biquadratics for the velocity, bilinears for the pressure, biquadratics for the elastic stress tensor, and bilinears for enforcing the rate of deformation (the first, second, third, and fourth equations in \eqref{eq:evss_system_D}, respectively). We are using a triangulated mesh, but we will try to follow their choices as closely as possible to discretize \eqref{eq:evss_system_D}: we choose continuous piecewise quadratic elements for the velocity, continuous piecewise linear elements for the pressure, continuous piecewise quadratic elements for the stress, and continuous piecewise linear elements for enforcing the definition of $\D$. Note that the velocity-pressure satisfy the usual LBB conditions. 

Letting $h$ denote a characteristic mesh size and subscript $h$ denote a function defined on a discrete mesh, the finite element formulation of \eqref{eq:evss_system_D}, with SUPG stabilization on the constitutive equation, says to find $(\vu_{h}, p_{h}, \Sig_{h}, \D_{h}) \in (P^{2}_{0}(\domain_{h})^{d} + \g)\times (P^{1}(\domain_{h})/\R) \times P^{2}(\domain_{h})^{d^{2}} \times P^{1}(\domain_{h})^{d^{2}}$, where $\g$ is the velocity Dirichlet data, such that 
\begin{equation}
    \begin{aligned}
        \int_{\domain}\Bigl(\eta_{0}\del\vu_{h}:\del\vv_{h} + (\vu_{h}\cdot\del\vu_{h})\cdot\vv_{h} +\del p_{h}  - \del\cdot\Sig_{h} - \f_{h}\Bigr)\cdot\vv_{h}\,dx &=0, \\
        \int_{\domain}(\del\cdot\vu_{h}) q_{h}\,dx &= 0, \\
        \int_{\domain} \bigg(\Sig_{h} + \lambda_{1}\Bigl(\vu_{h}\cdot\del(\Sig_{h} + 2\eta_{0}\D_{h}) - (\del\vu_{h})\cdot(\Sig_{h} + 2\eta_{0}\D_{h}) - (\Sig_{h} + 2\eta_{0}\D_{h})\cdot(\del\vu_{h})^{t}\Bigr) &\\
        +(\lambda_{1}-\mu_{1}) \Bigl(\D(\vu_{h})\cdot(\Sig_{h} + 2\eta_{0}\D_{h}) + (\Sig_{h} + 2\eta_{0}\D_{h})\cdot\D(\vu_{h})\Bigr) \bigg):\Bigl(\ess_{h} + h\vu_{h}\cdot\del\ess_{h}\Bigr) \,dx &= 0, \\
        \int_{\domain} \Bigl(2\D_{h}- \del\vu_{h} - (\del\vu_{h})^{t}\Bigr):\Phi_{h}\,dx &=0,
    \end{aligned}
    \label{eq:evss_system_D_fem}
\end{equation}
for all $(\vv_{h}$, $q_{h}$, $\ess_{h}, \Phi_{h})\in P^{2}_{0}(\domain_{h})^{d} \times P^{1}(\domain_{h}) \times P^{2}(\domain_{h})^{d^{2}} \times P^{1}(\domain_{h})^{d^{2}}$. Again, recall that, like in \eqref{eq:evss_system_D}, $\D_{h}$ is an independent variable to be solved for, while $\D(\vu_{h}) = (\del\vu_{h} + \del\vu_{h}^{t})/2$ is used explicitly. Since $\Sig$ is symmetric and $\D_{h}$ is symmetric and traceless, computation can be saved in practice by using smaller spaces than $P^{2}(\domain_{h})^{d^{2}}$ and $P^{1}(\domain_{h})^{d^{2}}$. We use Newton's method to address the nonlinearity of the entire mixed system \eqref{eq:evss_system_D_fem}. 

\subsection{SRTD}\label{subsec2:2_SRTD}
Renardy in 1985 was one of the first to give an existence proof for viscoelastic fluids, and he did so for the UCM model \cite{renardy_iter}. In his seminal paper, the divergence of the constitutive equation is taken and expanded out, and the divergence of $\T$ from the momentum equation is substituted in. A Picard-like iterative scheme is then described which alternates between solving a Stokes-like equation and the constitutive equation, and this effectively decouples the system. His iterative scheme also provided the basis for one of the first potential solution algorithms for a viscoelastic fluid. The scheme is now often referred to as \textit{Renardy iteration}. 

Girault and Scott in \cite{scott_wellposed} (summarized in \cite{scott_summary}) modify Renardy's method and generalize it to the O3 model. They too substitute the divergence of the stress tensor from the momentum equation into the divergence of the constitutive equation, but do so more selectively than Renardy. Hence, they coin their method ``Selective Replacement of Tensor Divergence," or SRTD (to be pronounced ``sorted" and not ``sordid"). Since deriving the SRTD formulation requires taking an extra derivative, a bulk of the paper \cite{scott_wellposed} is spent proving that the SRTD formulation is equivalent to the original O3 model given in \eqref{eq:oldroyd3_system}, assuming a connected Lipschitz domain $\domain$ of sufficient regularity (see \cite{scott_summary, scott_wellposed}). Specifically, \cite[Theorem 3.3]{scott_wellposed} says that for  of dimension $d\in\{2,3\}$, if $q>d$ and if $\vu\in W^{2}_{q}(\domain)^{d}$, $\T\in W^{1}_{q}(\domain)^{d^{2}}$, $f\in W^{1}_{q}(\domain)^{d}$, and $p\in W^{1}_{q}(\domain)/\R$, then the O3 forumulation \eqref{eq:oldroyd3_system} is equivalent to the SRTD formulation below:
\begin{equation}
    \begin{aligned}
        -\eta_{0}\Delta\vu + \vu\cdot\del\vu + \del\pi &= \F(\f,\vu,p,\T),  \\
        \del\cdot\vu &= 0, \\
        p + \lambda_{1}\vu\cdot\del p &= \pi, \\ 
        \T + \lambda_{1}( \vu\cdot\del\T + \W\cdot\T - \T\cdot\W) - \mu_{1} (\D\cdot\T + \T\cdot\D) &= 2\eta_{0}\D,
    \end{aligned}
    \label{eq:SRTD_formulation}
\end{equation}
where $\F$ is given by
\begin{equation}
    \begin{aligned}
        \F(\f, \vu, p, \T) = &\f + \lambda_{1}\vu\cdot\del \f + \lambda_{1}(\del\vu)^{t}\cdot\del p - \lambda_{1}\bigl(\vu\cdot\del(\vu\cdot\del\vu) - \del\cdot(\del\vu\cdot\T)\bigr) \\
        &- (\lambda_{1}-\mu_{1})\del\cdot(\D(\vu)\cdot\T + \T\cdot\D(\vu)),
    \end{aligned}
    \label{eq:F_SRTD}
\end{equation}
and $\pi$ is an auxiliary pressure function from which the true pressure $p$ can be recovered. Similar to Renardy iteration, Girault and Scott propose the following 3-stage iterative algorithm, effectively decoupling \eqref{eq:SRTD_formulation}. Given $\vu^{(n-1)}, p^{(n-1)}, \pi^{(n-1)}$, and $\T^{(n-1)}$, in stage 1, solve the Navier-Stokes--like system for $\vu^{(n)}$ and $\pi^{(n)}$: 
\begin{equation}
    \begin{aligned}
        -\eta_{0}\Delta\vu^{(n)} + \vu^{(n)}\cdot\del\vu^{(n)} + \del\pi^{(n)} &= \F(\f,\vu^{(n-1)},p^{(n-1)},\T^{(n-1)}),  \\
    \del\cdot\vu^{(n)} &= 0.
    \end{aligned}
    \label{eq:SRTD_stage1}
\end{equation}
In the second stage, given $\vu^{(n)}$ and $\pi^{(n)}$, solve the pressure transport equation for $p^{(n)}$:
\begin{equation}
    p^{(n)}+\lambda_{1} \vu^{(n)}\cdot \del p^{(n)} = \pi^{(n)}.
    \label{eq:SRTD_stage2}
\end{equation}
Finally, in stage 3, given $\vu^{(n)}$, solve the constitutive equation for $\T^{(n)}$: 
\begin{equation}
    \T^{(n)} + \lambda_{1}( \vu^{(n)}\cdot\del\T^{(n)} + \W^{(n)}\cdot\T^{(n)} - \T^{(n)}\cdot\W^{(n)}) - \mu_{1} (\D^{(n)}\cdot\T^{(n)} + \T^{(n)}\cdot\D^{(n)}) = 2\eta_{0}\D^{(n)},
    \label{eq:SRTD_stage3}
\end{equation}
where $\D^{(n)} = \D(\vu^{(n)})$ and $\W^{(n)}= \W(\vu^{(n)})$. In \cite[Section 6]{scott_wellposed}, the iterative scheme defined above for ($\vu^{(n)}, \pi^{(n)}, p^{(n)}, \T^{(n)})$ was shown to converge in $H^{1}(\domain)^{d}\times L_{2}(\domain)\times L_{2}(\domain)\times L_{2}(\domain)^{d^{2}}$ to the solution of \eqref{eq:SRTD_formulation} starting from $0$ initial guesses, $\vu^{(0)}=\mathbf{0}$, $p^{(0)}=0$, and $\T^{(0)}=\mathbf{0}$, and assuming sufficiently small non-Newtonian parameters $\lambda_{1}$ and $\mu_{1}$. The Weissenberg number for viscoelastic flow problems increases with the relaxation time $\lambda_{1}$, and as expected, the SRTD algorithm will fail to converge past certain Weissenberg numbers. 

The assumption of zero Dirichlet boundary conditions on the velocity was also assumed throughout \cite{scott_wellposed}, though the authors remark that this ``can be relaxed to allow $\vu=\g$ on $\partial\domain$ provided $\g\cdot \vn =0$, where $\vn$ is the outward normal to $\partial\domain$." The requirement of tangential boundary conditions will be discussed in the next section. 

\subsubsection{SRTD Finite Element Implementation}\label{subsubsec2:2_1_SRTD_FE}
Since the first stage of SRTD \eqref{eq:SRTD_stage1} requires a Navier-Stokes solve, we opted to use Taylor-Hood elements: piecewise quadratics for the velocity $\vu$ and piecewise linears for the auxiliary pressure $\pi$. This also matches the velocity-pressure pair we use to discretize the EVSS formulation. To get the weak form of Stage 1 of SRTD \eqref{eq:SRTD_stage1}, integration by parts is performed on $\mathcal{F}$ and $\Delta\vu$, and so the finite element formulation of \eqref{eq:SRTD_stage1} says to find $(\vu_{h}^{(n)}, \pi_{h}^{(n)})\in (P^{2}_{0}(\domain_{h})+\g)\times (P^{1}(\domain_{h})/\R)$ such that
\begin{equation}
    \begin{aligned}
        \eta_{0}\int_{\domain} \del\vu_{h}^{(n)}:\del\vv_{h}\,dx + &\int_{\domain}\Bigl (\vu_{h}^{(n)}\cdot\del\vu_{h}^{(n)} \Bigr)\cdot\vv_{h}\,dx - \int_{\domain}\pi_{h}^{(n)}\del\cdot\vv_{h}\,dx \\
        = &\int_{\domain} \f\cdot  \Bigl(\vv_{h}-\lambda_{1}\vu_{h}^{(n-1)}\cdot\del\vv_{h}\Bigr) \,dx - \lambda_{1}\bigg( \int_{\domain} p_{h}^{(n-1)}\Bigl(\del\vu_{h}^{(n-1)}\Bigr)^{t}:\del\vv_{h}\,dx \\
        &-\int_{\domain}\Bigl(\vu_{h}^{(n-1)}\cdot\del\vu_{h}^{(n-1)}\Bigr)\cdot \Bigl(\vu_{h}^{(n-1)}\cdot\del\vv_{h}\Bigr)\,dx \\
        & + \int_{\domain}\Bigl(\del\vu_{h}^{(n-1)}\cdot\T_{h}^{(n-1)}\Bigr):\del\vv_{h}\,dx \bigg) \\
        +&(\lambda_{1}-\mu_{1})\int_{\domain} \Bigl(\D_{h}^{(n-1)}\cdot\T_{h}^{(n-1)} + \T_{h}^{(n-1)}\cdot\D_{h}^{(n-1)}\Bigr):\del\vv_{h}\,dx,
    \end{aligned}
\end{equation}
where $\D_{h}^{(n-1)} = \D(\vu_{h}^{(n-1)})$, coupled with 
\begin{equation}
    \int_{\domain} \Bigl(\del\cdot\vu_{h}^{(n)}\Bigr)q_{h}\,dx = 0,
\end{equation}
for all $(\vv_{h}, q_{h})\in P^{2}_{0}(\D_{h})^{d}\times P^{1}(\domain_{h})$. 

Piecewise linears are used for the true pressure $p$ to solve \eqref{eq:SRTD_stage2}, matching the space we use in EVSS for the pressure. Stage 2 of SRTD says to find $p_{h}^{(n)}\in P^{1}(\domain_{h})$ such that 
\begin{equation}
    \int_{\domain}\Bigl (p^{(n)}_{h}+\lambda_{1}\vu^{(n)}_{h}\cdot\del p^{(n)}_{h}\Bigr) r_{h}\,dx = \int_{\domain}\pi^{(n)}_{h}r_{h}\,dx
    \label{eq:SRTD_stage2_fem}
\end{equation}
for all $r_{h}\in P^{1}(\domain_{h})$. 

Finally, matching our EVSS discretization, piecewise quadratics are used for the stress tensor $\T$ in \eqref{eq:SRTD_stage3}, and the finite element formulation of Stage 3 of SRTD says to find $\T^{(n)}_{h}\in P^{2}(\domain_{h})^{d^{2}}$ such that
\begin{align}
    \int_{\domain}\Big( \T_{h}^{(n)} &+ \lambda_{1}\Bigl( \vu_{h}^{(n)}\cdot\del\T_{h}^{(n)} + \W_{h}^{(n)}\cdot\T_{h}^{(n)} - \T_{h}^{(n)}\cdot\W_{h}^{(n)}\Bigr) \nonumber \\ 
    &- \mu_{1} \Bigl(\D_{h}^{(n)}\cdot\T_{h}^{(n)} + \T_{h}^{(n)}\cdot\D_{h}^{(n)}\Bigr) \Big):\ess_{h} \,dx = 2\eta_{0} \int_{\domain} \D_{h}^{(n)}:\ess_{h} \,dx
    \label{eq:SRTD_stage3_fem}
\end{align}
where $\W_{h}^{(n)} = \W(\vu_{h}^{(n)})$, for all $\ess\in P^{2}(\domain_{h})^{d^{2}}$. Like the EVSS finite element formulation, in practice, a smaller space than $P^{2}(\D_{h})^{d^{2}}$ can be used since $\T$ is symmetric. 

In preliminary testing, we ran SRTD with SUPG stabilization on the pressure transport equation \eqref{eq:SRTD_stage2_fem} and constitutive equation \eqref{eq:SRTD_stage3_fem}. In our test cases, we saw no notable change in performance in SRTD with SUPG versus without SUPG. This is to be expected, though, as we do not examine cases with particularly high Reynolds numbers, and SUPG could still be beneficial for advection dominated problems.

\section{Test Problems}\label{sec3:test_problems}
As mentioned, tangential Dirichlet boundary conditions were assumed in \cite{scott_wellposed}, and this requirement restricts the number of allowable flow problems by excluding those with nontrivial inflow and outflow. In this section, we describe two test problems which do satisfy the tangential Dirichlet data requirement: the eccentric rotating cylinders problem, sometimes called the journal-bearing problem, and the lid-driven cavity problem in both 2 and 3 dimensions. Due to limited computational resources and the large number of degrees of freedom that come with three-dimensional simulations, numerical results for the 3D lid-driven cavity problem are somewhat limited, but we do include some to show that SRTD is capable of handling three-dimensional problems. For all test problems, $U$ will denote a characteristic (scalar-valued) speed for the problem.

\subsection{The Journal-Bearing Problem}\label{subsec3:1_jb}
The eccentric rotating cylinders problem is characterized by two circles, the smaller of which has radius $r$ and is contained within the larger one of radius $R>r$, with centers offset from each other by an eccentricity value $e$. As a benchmark flow problem, typically the inner circle is set to rotate at constant angular speed, while the outer one is held fixed, and both are prescribed a no-slip boundary condition. This problem has application in the study of lubrication, where it is called the journal-bearing problem. In this context, the journal, represented by the inner circle, is free to spin within a bearing, represented by the larger circle, and the assembly is lubricated by a fluid; the fluid fills the gap between the two circles, preventing them from making direct contact, and is often described by a non-Newtonian model.

\begin{figure}[htbp!]
    \centering
    \begin{subfigure}{.5\textwidth}
        \centering
        \includegraphics[height=.84\linewidth]{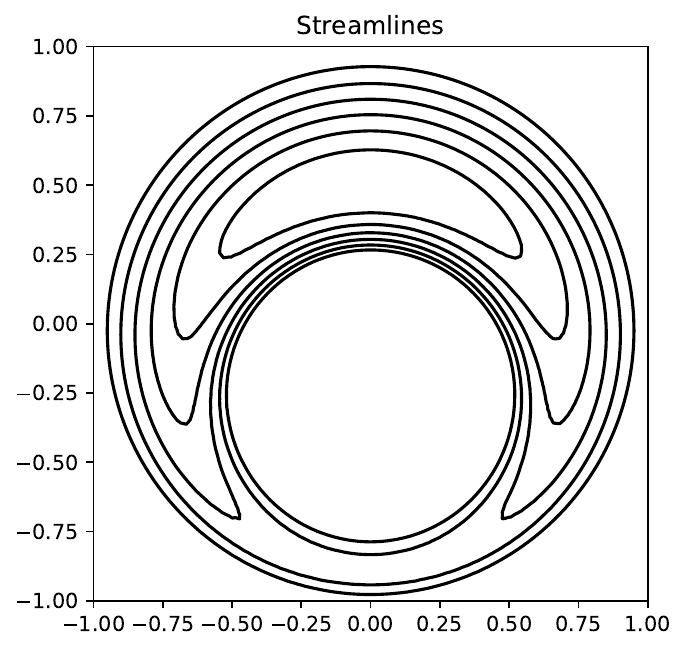}
    \end{subfigure}%
    \begin{subfigure}{.5\textwidth}
        \centering
        \includegraphics[height=.84\linewidth]{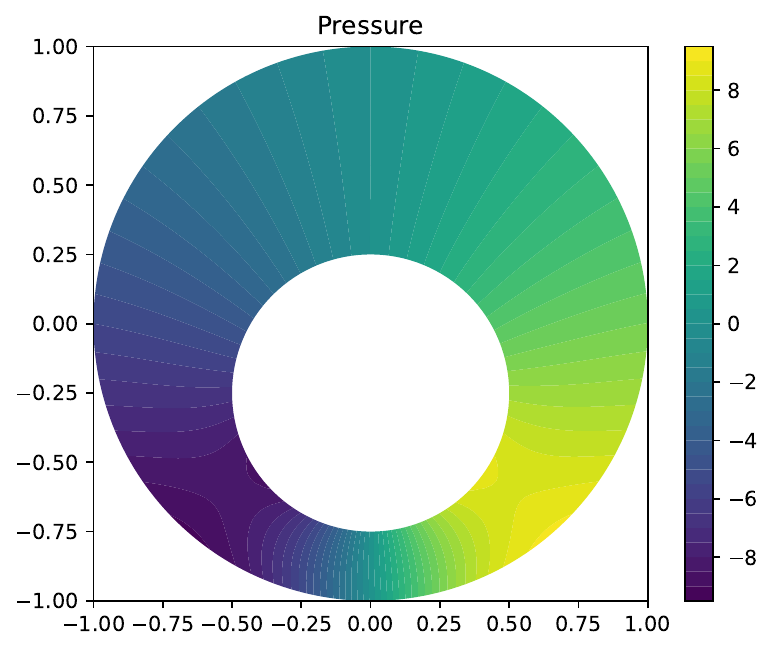}
    \end{subfigure}
    \caption{Streamlines (left) and the pressure profile (right) of the journal-bearing problem, with outer radius $R=1$, inner radius $r=0.5$, and eccentricity $e=0.25$, for the UCM model with $\lambda_{1}=0.01$, solved with the SRTD algorithm.}
    \label{fig:jb}
\end{figure}

The tangential speed of the inner circle, or journal, will serve as our characteristic speed $U$, giving the inner circle an angular velocity of $\omega=U/r$. In alignment with other literature (see, for instance \cite{vef_erc_huang_ptt,erc_ja_mechanics}), we define the Weissenberg number $\Wi$ and Reynolds number $\Rey$ for this problem as 
\begin{equation}
    \Wi := \lambda_{1}\omega = \lambda_{1}U/r, \quad \Rey := R^{2}\omega /\eta_{0} = R^{2}U/(r\eta_{0}).
    \label{wi_re_jb}
\end{equation}
We restrict to the case where the outer circle or bearing has radius $R=1$, the inner circle or journal has radius $r=0.5$, and the eccentricity is $e=0.25$. We will also restrict to the case where $\eta_{0}=1$. Thus, for the remainder of the paper, the Weissenberg number for the journal-bearing problem will be defined as $\Wi = 2\lambda_{1}U$, and the Reynolds number $\Rey = 2U$. An example of the streamlines and pressure profile for the journal-bearing problem, solved using the SRTD algorithm, can be seen in Figure \ref{fig:jb}. 

\subsection{The Lid-Driven Cavity Problem}\label{subsec3:2_ldc}
The lid-driven cavity problem is a popular test problem in computational fluid dynamics \cite{layton}. The geometry of the lid-driven cavity problem is given by a simple square in two dimensions or cube in three dimensions representing a physical cavity filled with fluid. The sides and the bottom floor are fixed, while the top wall slides in a certain direction and initiates a flow in the fluid. This can be described via a tangential Dirichlet boundary condition for the velocity on the top wall, and zero Dirichlet boundary conditions for the velocity on the bottom floor and the walls. The top wall can be set to move at constant speed, but this has the effect of introducing a discontinuity in the velocity profile at the top corners. In preliminary testing, this severely limited the performance of SRTD, which makes sense given the smoothness requirements for convergence. Instead, as is often done in the literature (see, for instance, \cite{ldc_of_vef_sousa}), we choose avoid this complication by prescribing a velocity profile along the top wall which decays as it reaches the corners. In the case where the cavity is the unit square $[0,1]^{2}$, the top lid Dirichlet data is given by
\[g_{top}(x,y) = (16x^{2}(1-x)^{2},\quad 0),\]
while in the three-dimensional unit cube case, we have
\[g_{top}(x,y,z) = (256x^{2}y^{2}(1-x)^{2}(1-y)^{2}, \quad 0, \quad 0).\] 
The coefficients $16$ and $256$ are chosen so the maximum velocity along the top wall is $1$, occurring at the midpoint in both cases. We further scale the top lid velocity up by a characteristic speed $U$, which in this context corresponds to the maximum velocity along the top lid. 

In the more general two-dimensional case where the cavity has height $H$ and width $L$, we define Weissenberg number $\Wi$ and the Reynolds number $\Rey$ as 
\begin{equation}
    \Wi = \lambda_{1}U/H, \quad \Rey=UH/\eta_{0}.
    \label{wi_re_ldc}
\end{equation}
However, we will only consider the unit square $[0,1]^{2}$ and unit cube $[0,1]^{3}$. So for the remainder of the paper, the Weissenberg number for the lid-driven cavity problem will be defined to be $\Wi = \lambda_{1}U$, and the Reynolds number $\Rey=U$. 

\begin{figure}[htbp!]
    \centering
    \begin{subfigure}{.5\textwidth}
        \centering
        \includegraphics[height=.84\linewidth]{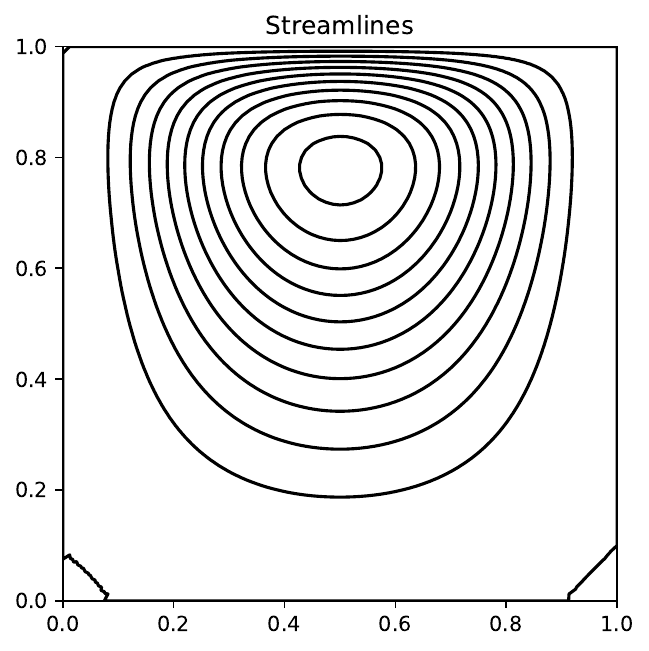}
    \end{subfigure}%
    \begin{subfigure}{.5\textwidth}
        \centering
        \includegraphics[height=.84\linewidth]{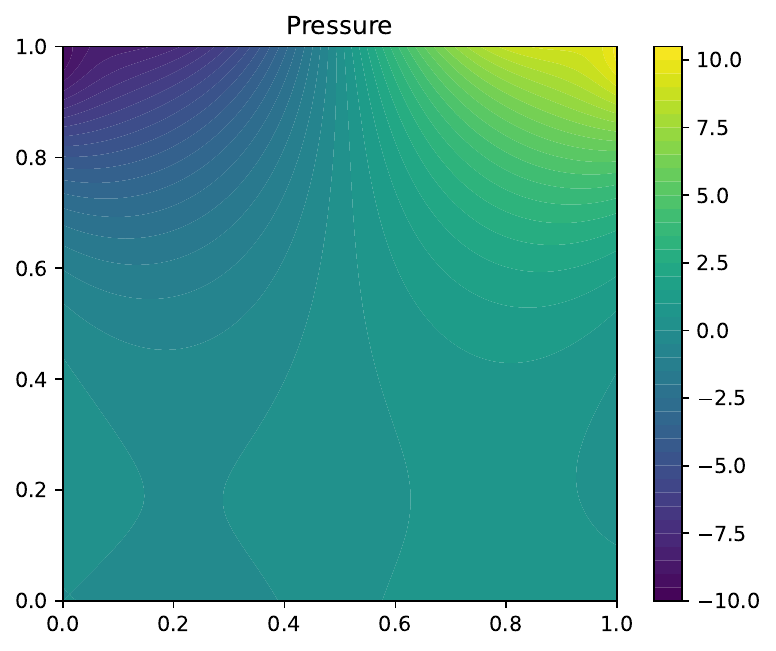}
    \end{subfigure}
    \caption{Streamlines (left) and the pressure profile (right) of the lid-driven cavity problem for the UCM model with $\lambda_{1}=0.01$.}
    \label{fig:ldc}
\end{figure}

An example of the streamlines and pressure profile of the two-dimensional lid-driven cavity problem, solved using the SRTD algorithm, can be seen in Figure \ref{fig:ldc}.

\section{Numerical Results}\label{sec4:numerical_results}
In this section, we give some numerical results, beginning with the approximate Weissenberg number limits for both SRTD and EVSS across both the journal-bearing problem and lid-driven cavity problem and with both the UCM and corotational Maxwell models. We further examine how the SRTD algorithm converges on a per-iteration basis for various $\lambda_{1}$ values, corresponding to increasing Weissenberg number. In the subsections after, we examine how SRTD performs with regards to mesh refinement. Since, as we will see, the SRTD algorithm performs well only for small Weissenberg numbers, we close the section by examining how much the non-Newtonian flow differs from Newtonian flow for $\lambda_{1}$ values for which SRTD does converge. 

\begin{figure}[htbp!]
    \centering
    \begin{subfigure}{.5\textwidth}
        \centering
        \includegraphics[height=.84\linewidth]{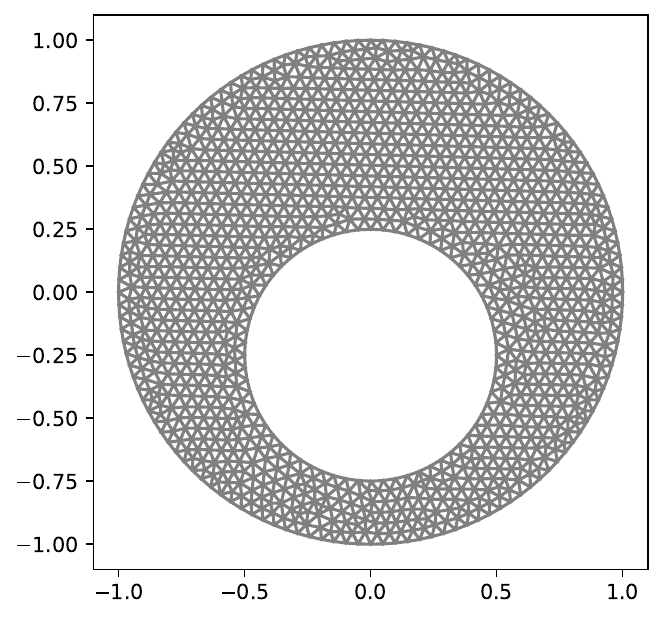}
    \end{subfigure}%
    \begin{subfigure}{.5\textwidth}
        \centering
        \includegraphics[height=.84\linewidth]{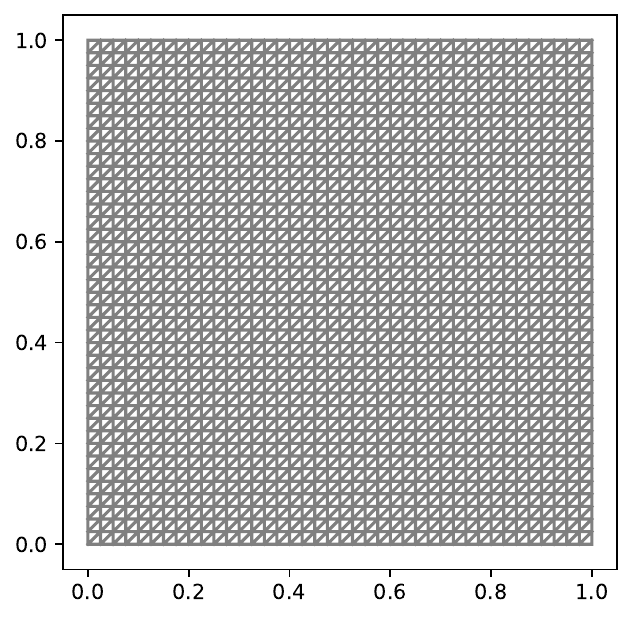}
    \end{subfigure}
    \caption{Gmsh-generated mesh for the journal-bearing problem with characteristic mesh size $h=0.05$ (left) and the unit square mesh generated by FEniCS with mesh size parameter $40$, or characteristic mesh size $h=1/40 = 0.025$ (right).}
    \label{fig:meshes}
\end{figure}

We obtained numerical results via the Python finite element system FEniCS \cite{fenics}. The meshes for the lid-driven cavity problem in two and three dimensions were generated using FEniCS' built-in \verb|UnitSquareMesh()| and \verb|UnitCubeMesh()| functions, respectively. The meshes for the journal-bearing problem were generated using Gmsh \cite{gmsh}. An example of a Gmsh-generated mesh for the journal-bearing problem with characteristic mesh size $h=0.05$ and the FEniCS' generated unit square mesh for the lid-driven cavity problem with characteristic mesh size $h=0.025$ can be found in Figure \ref{fig:meshes}. The number of degrees of freedom for the 3 SRTD stages and the full EVSS solve for various characteristic mesh sizes are listed in Tables \ref{tab:jb_dof}, \ref{tab:ldc2d_dof}, and \ref{tab:ldc3d_dof} for the journal-bearing problem, the 2D lid-driven cavity problem, and the 3D lid-driven cavity problem, respectively, when using the finite element spaces described in Section 2. 

\begin{table}[htbp!]
\centering
\begin{tabular}{rr|rrr|r}
\multicolumn{2}{r|}{\textbf{Journal-Bearing}}                 & \multicolumn{3}{c|}{\textbf{SRTD}}                                                & \multicolumn{1}{c}{\textbf{EVSS}} \\ \hline
\multicolumn{1}{r|}{$h$}        & \multicolumn{1}{c|}{Elements} & \multicolumn{1}{r|}{Stage 1 DoF} & \multicolumn{1}{r|}{Stage 2 DoF} & Stage 3 DoF & Total DoF                         \\ \hline
\multicolumn{1}{r|}{2.00e-01} & 162                           & \multicolumn{1}{r|}{657}         & \multicolumn{1}{r|}{105}         & 1116        & 2175                              \\
\multicolumn{1}{r|}{1.00e-01} & 606                           & \multicolumn{1}{r|}{2583}        & \multicolumn{1}{r|}{351}         & 3924        & 7593                              \\
\multicolumn{1}{r|}{5.00e-02} & 2350                          & \multicolumn{1}{r|}{10287}       & \multicolumn{1}{r|}{1271}        & 14676       & 28273                             \\
\multicolumn{1}{r|}{2.50e-02} & 9000                          & \multicolumn{1}{r|}{39930}       & \multicolumn{1}{r|}{4690}        & 55140       & 105970                            \\
\multicolumn{1}{r|}{1.25e-02} & 35328                         & \multicolumn{1}{r|}{157842}      & \multicolumn{1}{r|}{18042}       & 214236      & 411186                            \\ \hline
\end{tabular}
\caption{Total number of elements (triangles) for the journal-bearing problem for various mesh sizes, plus the number of total degrees of freedom for Stage 1 (Navier-Stokes equations for velocity and auxiliary pressure), Stage 2 (pressure transport equation for the true pressure), and Stage 3 (constitutive equation for the stress tensor) of SRTD, plus the total number of degrees of freedom for the full mixed EVSS formulation. Recall that Stage 1 of SRTD and the full EVSS formulation are nonlinear and require Newton's method to be solved.}
\label{tab:jb_dof}
\end{table}

\begin{table}[htbp!]
\centering
\begin{tabular}{rr|rrr|r}
\multicolumn{2}{r|}{\textbf{2D Lid-Driven Cavity}} & \multicolumn{3}{c|}{\textbf{SRTD}}                                                & \multicolumn{1}{c}{\textbf{EVSS}} \\ \hline
\multicolumn{1}{r|}{$h$}             & Elements      & \multicolumn{1}{r|}{Stage 1 DoF} & \multicolumn{1}{r|}{Stage 2 DoF} & Stage 3 DoF & Total DoF                         \\ \hline
\multicolumn{1}{r|}{1.00e-01}      & 200           & \multicolumn{1}{r|}{1003}        & \multicolumn{1}{r|}{121}         & 1323        & 2568                              \\
\multicolumn{1}{r|}{5.00e-02}      & 800           & \multicolumn{1}{r|}{3803}        & \multicolumn{1}{r|}{441}         & 5043        & 9728                              \\
\multicolumn{1}{r|}{2.50e-02}      & 3200          & \multicolumn{1}{r|}{14803}       & \multicolumn{1}{r|}{1681}        & 19683       & 37848                             \\
\multicolumn{1}{r|}{1.25e-02}      & 12800         & \multicolumn{1}{r|}{58403}       & \multicolumn{1}{r|}{6561}        & 77763       & 149288                            \\
\multicolumn{1}{r|}{6.25e-03}      & 51200         & \multicolumn{1}{r|}{232003}      & \multicolumn{1}{r|}{25921}       & 309123      & 592968                            \\ \hline
\end{tabular}
\caption{Total number of elements (triangles) for the two-dimensional lid-driven cavity problem for various mesh sizes, plus the number of total degrees of freedom for Stage 1 (Navier-Stokes equations for velocity and auxiliary pressure), Stage 2 (pressure transport equation for the true pressure), and Stage 3 (constitutive equation for the stress tensor) of SRTD, and the total number of degrees of freedom for the full mixed EVSS formulation. Recall that Stage 1 of SRTD and the full EVSS formulation are nonlinear and require Newton's method to be solved.}
\label{tab:ldc2d_dof}
\end{table}

\begin{table}[htbp!]
\centering
\begin{tabular}{rr|rrr|r}
\multicolumn{2}{r|}{\textbf{3D Lid-Driven Cavity}} & \multicolumn{3}{c|}{\textbf{SRTD}}                                                & \multicolumn{1}{c}{\textbf{EVSS}} \\ \hline
\multicolumn{1}{r|}{$h$}             & Elements      & \multicolumn{1}{r|}{Stage 1 DoF} & \multicolumn{1}{r|}{Stage 2 DoF} & Stage 3 DoF & Total DoF                         \\ \hline
\multicolumn{1}{r|}{2.50e-01}      & 384           & \multicolumn{1}{r|}{2312}        & \multicolumn{1}{r|}{125}         & 4374        & 7311                              \\
\multicolumn{1}{r|}{1.25e-01}      & 3072          & \multicolumn{1}{r|}{15468}       & \multicolumn{1}{r|}{729}         & 29478       & 48591                             \\
\multicolumn{1}{r|}{6.25e-02}      & 24576         & \multicolumn{1}{r|}{112724}      & \multicolumn{1}{r|}{4913}        & 215622      & 352911                            \\ \hline
\end{tabular}
\caption{Total number of elements (tetrahedra) for the three-dimensional lid-driven cavity problem for various mesh sizes, plus the number of total degrees of freedom for Stage 1 (Navier-Stokes equations for velocity and auxiliary pressure), Stage 2 (pressure transport equation for the true pressure), and Stage 3 (constitutive equation for the stress tensor) of SRTD, and the total number of degrees of freedom for the full mixed EVSS formulation. Recall that Stage 1 of SRTD and the full EVSS formulation are nonlinear and require Newton's method to be solved. }
\label{tab:ldc3d_dof}
\end{table}

The Stage 1 of the SRTD algorithm and the entire mixed EVSS formulation were solved using standard Newton's method as it is implemented in \verb|NonlinearVariationalSolver()| in FEniCS. A relative tolerance of $10^{-9}$ was used for the Newton iterations. As a starting guess for EVSS' Newton iteration, we used a Navier-Stokes solution velocity and pressure, but zero for the elastic stress tensor and deformation. The SRTD iteration officially starts from a zero initial guess, but as one can observe from \eqref{eq:SRTD_formulation} and \eqref{eq:F_SRTD}, when there are no body forces, then the first stage of the first iteration of SRTD \eqref{eq:SRTD_stage1} produces a Navier-Stokes solution with the same zero body forces. The previous SRTD iteration velocity $\vu^{(n-1)}$ and auxiliary pressure $\pi^{(n-1)}$ are used as the starting guess for the Newton iteration required by Stage 1 of the next SRTD iteration.

\subsection{Weissenberg Number Limits}\label{subsec4:1_Weissenberg}

We begin first by stating the observed approximate Weissenberg number limits for both SRTD and EVSS across both problems and for both the UCM model and the corotational Maxwell model, which are given in Table \ref{tab:weissenberg} below. The EVSS method is able to reach a significantly higher Weissenberg number for the UCM model on both problems than SRTD is able to, but only slightly higher for the corotational model. 

\begin{table}[htbp!]
\centering
\resizebox{0.4\textwidth}{!}{%
\begin{tabular}{c|cc|cc|}
     & \multicolumn{2}{c|}{UCM Model}   & \multicolumn{2}{c|}{Corotational} \\ \cline{2-5} 
     & \multicolumn{1}{c|}{LDC}  & JB   & \multicolumn{1}{c|}{LDC}       & JB       \\ \hline
SRTD & \multicolumn{1}{c|}{0.05} & 0.12 & \multicolumn{1}{c|}{0.06}      & 0.11     \\
EVSS & \multicolumn{1}{c|}{0.45} & 2.00 & \multicolumn{1}{c|}{0.09}      & 0.22     \\ \hline
\end{tabular}%
}
\caption{Approximate Weissenberg number limits for which the SRTD algorithm would converge and Newton iteration on EVSS would converge. }
\label{tab:weissenberg}
\end{table}

\subsection{SRTD Iteration Behavior}\label{subsec4:2_SRTD_iterates}

\begin{figure}[htbp!]
    \centering
    \begin{subfigure}{.5\textwidth}
        \centering
        \includegraphics[width=\linewidth]{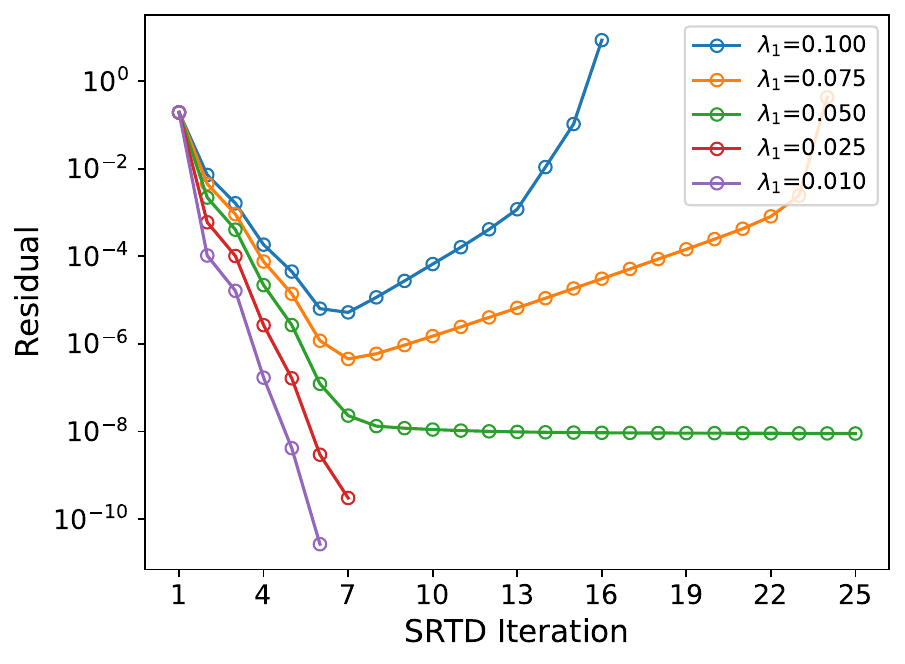}
    \end{subfigure}%
    \begin{subfigure}{.5\textwidth}
        \centering
        \includegraphics[width=\linewidth]{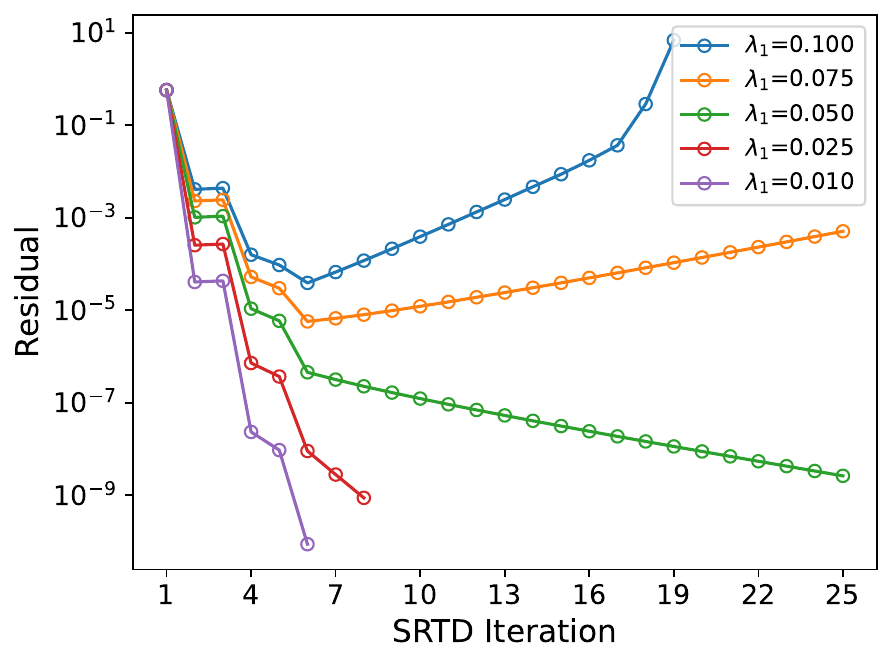}
    \end{subfigure}
    \caption{Number of SRTD iterations plotted against the $L_{2}$ difference between consecutive iterates for the UCM model ($\mu_{1}=\lambda_{1}$) with characteristic velocity $U=1$ on the two-dimensional lid-driven cavity problem (left) with mesh size $h=0.0125$ and the journal-bearing problem (right) with $h=0.025$.}
    \label{fig:ucm_num_iters_compare}
\end{figure}

\begin{figure}[htbp!]
    \centering
    \begin{subfigure}{.5\textwidth}
        \centering
        \includegraphics[width=\linewidth]{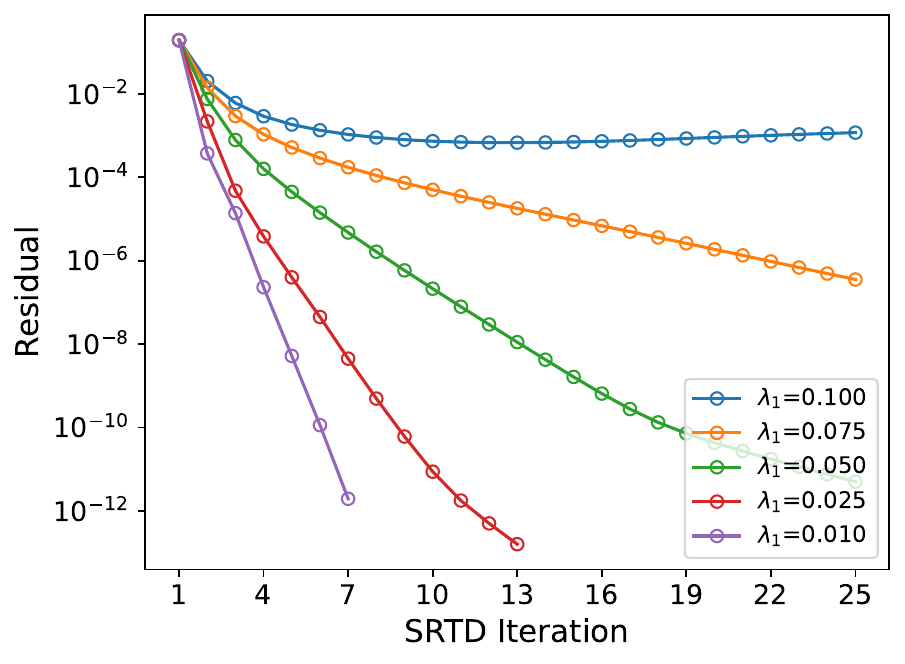}
    \end{subfigure}%
    \begin{subfigure}{.5\textwidth}
        \centering
        \includegraphics[width=\linewidth]{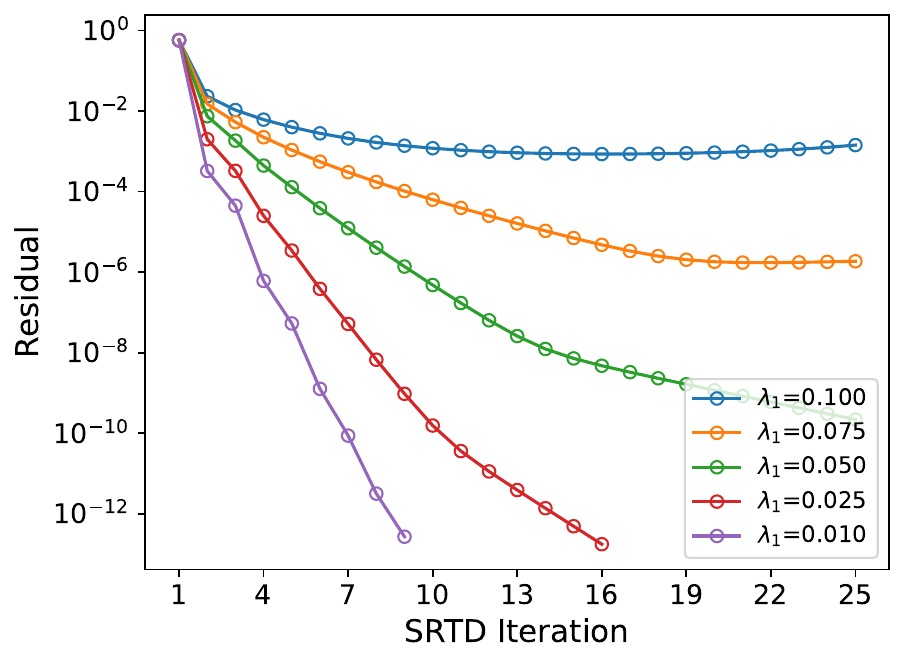}
    \end{subfigure}
    \caption{Number of SRTD iterations plotted against the $L_{2}$ difference between consecutive iterates for the corotational Maxwell model ($\mu_{1}=0$) with characteristic velocity $U=1$ on the lid-driven cavity problem (left) with mesh size $h=0.0125$ and the journal-bearing problem (right) with $h=0.025$.}
    \label{fig:corot_num_iters_compare}
\end{figure}

We say the SRTD algorithm converges if the relative error reaches a desired tolerance (we usually set it at $10^{-9}$) within the maximum number of SRTD iterations allowed (we usually used about 20). The Weissenberg number limit listed for SRTD in Table \ref{tab:weissenberg} does not appear to be a sharp cutoff. Indeed, for Weissenberg numbers near the limit listed in Table \ref{tab:weissenberg}, the relative error of the SRTD iterates will asymptote. For numbers significantly above the limit listed in Table \ref{tab:weissenberg}, the SRTD iterates may begin to decrease, but eventually increase. This behavior continues until the maximum number of SRTD iterations is reached, or the Newton iteration in the Navier-Stokes stage diverges. Figure \ref{fig:ucm_num_iters_compare} and Figure \ref{fig:corot_num_iters_compare} plot SRTD iteration's relative error against the iteration count for various $\lambda_{1}$ values near the Weissenberg number limit given in Table \ref{tab:weissenberg}. As one might expect, SRTD performance varies with the relaxation time, with performance largely decreasing as $\lambda_{1}$ increases.

\subsection{Mesh Refinement Behavior}\label{subsec4:3_mesh_ref}
With a sense for how SRTD behaves on fixed meshes, we would like to investigate its behavior as the mesh is refined. We also include mesh refinement experiments with EVSS for comparison. We consider the convergence of velocity profile $\vu$ in the $L_{2}$ and $H^{1}$ norms, and the $L_{2}$ convergence of the pressure $p$. Since the exact solution is unknown, we approximate the error by computing the norm of the difference between solutions on successively finer meshes, 
\begin{equation}
    \text{error}_{h/2} = \|\vu_{h} - \vu_{h/2}\|,
    \label{eq:approx_error}
\end{equation}
where $\vu_{h}$ denotes the calculated finite element solution on a mesh with characteristic mesh size $h$, and likewise for the pressure $p_{h}$. In FEniCS, this is computed by first interpolating the solution found on the coarser mesh $\vu_{h}$ onto the finer mesh using the \verb|interpolate()| function, and then computing the appropriate $L_{2}$ or $H^{1}$ norm using FEniCS' \verb|errornorm()|. The convergence rate for each method is then approximated using the following formula,
\begin{equation}
    \text{rate}_{h/4} = \log_{2}\left(\frac{\|\vu_{h} - \vu_{h/2}\|}{\|\vu_{h/2} - \vu_{h/4}\|}\right). 
    \label{eq:approx_convergence_rate}
\end{equation}
Since \eqref{eq:approx_error} requires a solution on two distinct mesh sizes, to obtain the data in Tables \ref{tab:ldc_srtd_ucm}-\ref{tab:ldc3d_srtd}, there is a coarsest mesh on which we must compute a solution but from which we can infer no error data, and their respective rows in the tables are omitted. Likewise, since \eqref{eq:approx_convergence_rate} requires three mesh sizes, for the second-coarsest mesh on which we compute a solution, we can compute an approximate error, but cannot compute an approximate rate. We include their respective rows row in Tables \ref{tab:ldc_srtd_ucm}-\ref{tab:ldc3d_srtd}, and the columns containing the approximate rates are left blank. 

We set the SRTD algorithm to terminate once a residual value of $10^{-9}$ was reached or after $20$ SRTD iterations. If the residual threshold was not reached within $20$ iterates, the algorithm would return the solution when the residual was the lowest.

\subsubsection{Convergence Rates for the Two-Dimensional Lid-Driven Cavity Problem}\label{subsubsec4:3_1_ldc_convergence}

When using degree $k$ polynomials, standard theory expects a convergence rate of $O(h^{k})$ in the $H^{1}$ norm and $O(h^{k+1})$ in the $L_{2}$ norm. This means that, since we implemented both SRTD and EVSS to use piecewise quadratics for the velocity and piecewise linears for the pressure, we should expect an $L_{2}$ convergence rate of $3$ and an $H^{1}$ convergence rate of $2$ for the velocity, and an $L_{2}$ convergence rate of $2$ for the pressure. For Weissenberg numbers smaller than those listed in Table \ref{tab:weissenberg}, we indeed see the expected convergence rates for both SRTD and EVSS. Some results for SRTD on the UCM model with $\Wi=0.05$ and corotational model with $\Wi=0.04$ can be found in Tables \ref{tab:ldc_srtd_ucm} and \ref{tab:ldc_srtd_corot}, respectively. Results for EVSS on the UCM model with $\Wi = 0.1$ and corotational model with $\Wi=0.04$ can be found in Tables \ref{tab:ldc_evss_ucm} and \ref{tab:ldc_evss_corot}, respectively. As we approach the approximate Weissenberg number limits, while the SRTD iterates and the Newton iteration on EVSS may still converge for a time, meaning they eventually reach a relative tolerance of $10^{-9}$, the methods do not achieve optimal convergence rates.

\begin{table}[htbp!]
\centering
\begin{tabular}{c|c|c|c|c|c|c}
\hline
$h$         & $L^{2}$ error $\vu$ & $L^{2}$ rate $\vu$ & $H^{1}$ error $\vu$ & $H^{1}$ rate $\vu$ & $L^{2}$ error $p$ & $L^{2}$ rate $p$ \\ \hline
5.00e-02 & 8.796e-04           & -                  & 6.474e-02           & -                  & 4.248e-02         & -                \\
2.50e-02 & 1.140e-04           & 2.948              & 1.710e-02           & 1.921              & 1.078e-02         & 1.979            \\
1.25e-02 & 1.415e-05           & 3.009              & 4.312e-03           & 1.987              & 2.614e-03         & 2.043            \\
6.25e-03 & 1.761e-06           & 3.006              & 1.078e-03           & 2.000              & 6.439e-04         & 2.022            \\ \hline
\end{tabular}
\caption{Mesh refinement for the SRTD formulation on the two-dimensional lid-driven cavity problem with $\Wi=0.05$ (with $U=1$ and $\lambda_{1}=0.05$), solved with the SRTD iteration, with $\eta_{0}=1$, and $\mu_{1}=\lambda_{1}$ (corresponding to the UCM model).}
\label{tab:ldc_srtd_ucm}
\end{table}

\begin{table}[htbp!]
\centering
\begin{tabular}{r|c|c|c|c|c|c}
\hline
$h$       & $L^{2}$ error $\vu$ & $L^{2}$ rate $\vu$ & $H^{1}$ error $\vu$ & $H^{1}$ rate $\vu$ & $L^{2}$ error $p$ & $L^{2}$ rate $p$ \\ \hline
5.00e-02 & 9.462e-04           & -                 & 6.746e-02           & -                  & 4.246e-02         & -                \\
2.50e-02 & 1.256e-04           & 2.913              & 1.828e-02           & 1.884              & 1.082e-02         & 1.973            \\
1.25e-02 & 1.592e-05           & 2.979              & 4.696e-03           & 1.961              & 2.680e-03         & 2.013            \\
6.25e-03 & 1.998e-06           & 2.995              & 1.183e-03           & 1.989              & 6.629e-04         & 2.015            \\ \hline
\end{tabular}
\caption{Mesh refinement for the SRTD formulation on the two-dimensional lid-driven cavity problem with $\Wi=0.04$ (with $U=1$ and $\lambda_{1}=0.04$), solved with the SRTD iteration, with $\eta_{0}=1$, and $\mu_{1}=0$ (corresponding to the corotational Maxwell model). }
\label{tab:ldc_srtd_corot}
\end{table}

\begin{table}[htbp!]
\centering
\begin{tabular}{c|c|c|c|c|c|c}
\hline
$h$         & $L^{2}$ error $\vu$ & $L^{2}$ rate $\vu$ & $H^{1}$ error $\vu$ & $H^{1}$ rate $\vu$ & $L^{2}$ error $p$ & $L^{2}$ rate $p$ \\ \hline
5.00e-02 & 1.224e-03           & -                  & 7.891e-02           & -                  & 6.232e-02         & -                \\
2.50e-02 & 1.598e-04           & 2.937              & 2.129e-02           & 1.890              & 1.552e-02         & 2.005            \\
1.25e-02 & 1.904e-05           & 3.069              & 5.352e-03           & 1.992              & 3.687e-03         & 2.074            \\
6.25e-03 & 2.304e-06           & 3.047              & 1.332e-03           & 2.007              & 9.002e-04         & 2.034            \\ \hline
\end{tabular}
\caption{Mesh refinement for the EVSS formulation on the two-dimensional lid-driven cavity problem with $\Wi=0.1$ (with $U=1$ and $\lambda_{1}=0.1$), solved fully coupled using Newton's method, with $\eta_{0}=1$, and $\mu_{1}=\lambda_{1}$ (corresponding to the UCM model).}
\label{tab:ldc_evss_ucm}
\end{table}

\begin{table}[htbp!]
\centering
\begin{tabular}{r|c|c|c|c|c|c}
\hline
$h$       & $L^{2}$ error $\vu$ & $L^{2}$ rate $\vu$ & $H^{1}$ error $\vu$ & $H^{1}$ rate $\vu$ & $L^{2}$ error $p$ & $L^{2}$ rate $p$ \\ \hline
5.00e-02 & 1.009e-03           & -                  & 6.998e-02           & -                  & 4.246e-02         & -                \\
2.50e-02 & 1.301e-04           & 2.955              & 1.863e-02           & 1.909              & 1.083e-02         & 1.971            \\
1.25e-02 & 1.616e-05           & 3.009              & 4.736e-03           & 1.976              & 2.683e-03         & 2.013            \\
6.25e-03 & 2.009e-06           & 3.008              & 1.188e-03           & 1.996              & 6.631e-04         & 2.016  \\
\hline      
\end{tabular}
\caption{Mesh refinement for the EVSS formulation on the two-dimensional lid-driven cavity problem with $\Wi=0.04$ (with $U=1$ and $\lambda_{1}=0.04$), solved fully coupled using Newton's method, with $\eta_{0}=1$, and $\mu_{1}=0$ (corresponding to the corotational Maxwell model).}
\label{tab:ldc_evss_corot}
\end{table}

\pagebreak 

\subsubsection{Convergence Rates for the Journal-Bearing Problem}\label{subsubsec4:3_2_jb_convergence}

The domain of the journal-bearing problem in the continuous case contains two circular boundaries, yet we are discretizing it with triangles. As described in \cite{scott_polygonal}, when a curved domain is approximated with polygons and quadratic piecewise polynomials are used to approximate the solution function, one expects an $H^{1}$ error bound of $O(h^{3/2})$ rather than the usual $O(h^{2})$, and $3/2$ is much closer to what is seen numerically in Tables \ref{tab:jb_srtd_ucm}, \ref{tab:jb_srtd_corot}, \ref{tab:jb_evss_ucm}, and \ref{tab:jb_evss_corot}. The pressure, which is approximated by piecewise linears, experiences the optimal and expected $L_{2}$ convergence rate of $2$. We recommend \cite{scott_polygonal} and Chapter 22 of \cite{scott_automated_modeling_fenics}, and the references therein, for more information on approximating curved domains with polygonal elements.

\begin{table}[htbp!]
\centering
\begin{tabular}{r|c|c|c|c|c|c}
\hline
$h$      & $L^{2}$ error $\vu$ & $L^{2}$ rate $\vu$ & $H^{1}$ error $\vu$ & $H^{1}$ rate $\vu$ & $L^{2}$ error $p$ & $L^{2}$ rate $p$ \\ \hline
1.00e-01 & 1.870e-02           & -                  & 2.746e-01           & -                  & 1.962e-01         & -                \\
5.00e-02 & 4.841e-03           & 1.950              & 9.908e-02           & 1.471              & 5.670e-02         & 1.791            \\
2.50e-02 & 1.220e-03           & 1.989              & 3.493e-02           & 1.504              & 1.515e-02         & 1.904            \\
1.25e-02 & 3.046e-04           & 2.002              & 1.289e-02           & 1.438              & 3.753e-03         & 2.013    \\
\hline
\end{tabular}
\caption{Mesh refinement for the SRTD formulation on the journal-bearing problem with $\Wi=0.1$ (with $U=1$ and $\lambda_{1}=0.05$), solved with the SRTD iteration, with $\eta_{0}=1$, and $\mu_{1}=\lambda_{1}$ (corresponding to the UCM model).}
\label{tab:jb_srtd_ucm}
\end{table}

\begin{table}[htbp!]
\centering
\begin{tabular}{r|c|c|c|c|c|c}
\hline
$h$      & $L^{2}$ error $\vu$ & $L^{2}$ rate $\vu$ & $H^{1}$ error $\vu$ & $H^{1}$ rate $\vu$ & $L^{2}$ error $p$ & $L^{2}$ rate $p$ \\ \hline
1.00e-01 & 1.834e-02           & -                  & 2.955e-01           & -                  & 1.703e-01         & -                \\
5.00e-02 & 4.798e-03           & 1.935              & 1.071e-01           & 1.465              & 5.061e-02         & 1.750            \\
2.50e-02 & 1.220e-03           & 1.976              & 3.755e-02           & 1.511              & 1.377e-02         & 1.878            \\
1.25e-02 & 3.055e-04           & 1.997              & 1.370e-02           & 1.454              & 3.474e-03         & 1.987      \\
\hline
\end{tabular}
\caption{Mesh refinement for the SRTD formulation on the journal-bearing problem with $\Wi=0.1$ (with $U=1$ and $\lambda_{1}=0.05$), solved with the SRTD iteration, with $\eta_{0}=1$, and $\mu_{1}=0$ (corresponding to the corotational Maxwell model).}
\label{tab:jb_srtd_corot}
\end{table}

\begin{table}[htbp!]
\centering
\begin{tabular}{r|c|c|c|c|c|c}
\hline
$h$       & $L^{2}$ error $\vu$ & $L^{2}$ rate $\vu$ & $H^{1}$ error $\vu$ & $H^{1}$ rate $\vu$ & $L^{2}$ error $p$ & $L^{2}$ rate $p$ \\ \hline
1.00e-01 & 3.159e-02           & -                  & 4.556e-01           & -                  & 7.305e-01         & -                \\
5.00e-02 & 5.614e-03           & 2.492              & 1.361e-01           & 1.743              & 1.126e-01         & 2.698            \\
2.50e-02 & 1.251e-03           & 2.166              & 4.269e-02           & 1.673              & 2.741e-02         & 2.038            \\
1.25e-02 & 3.086e-04           & 2.019              & 1.447e-02           & 1.561              & 6.545e-03         & 2.066  \\
\hline
\end{tabular}
\caption{Mesh refinement for the EVSS formulation on the journal-bearing problem problem with $\Wi=1.0$ (with $U=1$ and $\lambda_{1}=0.5$), solved fully coupled using Newton's method, with $\eta_{0}=1$, and $\mu_{1}=\lambda_{1}$ (corresponding to the UCM model).}
\label{tab:jb_evss_ucm}
\end{table}

\begin{table}[htbp!]
\centering
\begin{tabular}{r|c|c|c|c|c|c}
\hline
$h$       & $L^{2}$ error $\vu$ & $L^{2}$ rate $\vu$ & $H^{1}$ error $\vu$ & $H^{1}$ rate $\vu$ & $L^{2}$ error $p$ & $L^{2}$ rate $p$ \\ \hline
1.00e-01 & 1.805e-02           & -                  & 3.011e-01           & -                  & 1.657e-01         & -                \\
5.00e-02 & 4.739e-03           & 1.929              & 1.089e-01           & 1.467              & 4.937e-02         & 1.746            \\
2.50e-02 & 1.212e-03           & 1.967              & 3.784e-02           & 1.525              & 1.351e-02         & 1.869            \\
1.25e-02 & 3.044e-04           & 1.993              & 1.374e-02           & 1.462              & 3.408e-03         & 1.988 \\
\hline
\end{tabular}
\caption{Mesh refinement for the EVSS formulation on the journal-bearing problem problem with $\Wi=0.1$ (with $U=1$ and $\lambda_{1}=0.05$), solved fully coupled using Newton's method, with $\eta_{0}=1$, and $\mu_{1}=0$ (corresponding to the corotational Maxwell model).}
\label{tab:jb_evss_corot}
\end{table}

\subsubsection{Three-Dimensional Lid-Driven Cavity}\label{subsubsec4:3_3_ldc3d_convergence}
Due to the large number of degrees of freedom as the mesh is refined in three dimensions, we have only limited results for how SRTD performs here. Despite starting on a rather coarse mesh ($h=0.25$), the approximated convergence rates after only $2$ refinements are close to the expected optimal numbers, as one can see in Table \ref{tab:ldc3d_srtd}. This shows SRTD's viability for three-dimensional problems. 

\begin{table}[htbp!]
\centering
\begin{tabular}{r|c|c|c|c|c|c}
\hline
$h$       & $L^{2}$ error $\vu$ & $L^{2}$ rate $\vu$ & $H^{1}$ error $\vu$ & $H^{1}$ rate $\vu$ & $L^{2}$ error $p$ & $L^{2}$ rate $p$ \\ \hline
1.25e-01 & 9.287e-03           & -                  & 2.922e-01           & -                  & 2.905e-01         & -                \\
6.25e-02 & 1.427e-03           & 2.703              & 8.648e-02           & 1.757              & 5.523e-02         & 2.395   \\
\hline
\end{tabular}
\caption{Mesh refinement for SRTD formulation on the three-dimensional lid-driven cavity problem with $\Wi=0.01$ (with $U=1$ and $\lambda_{1}=0.01$), solved using the SRTD iteration, with $\eta_{0}=1$ and $\mu_{1}=\lambda_{1}$ (corresponding to the UCM model).}
\label{tab:ldc3d_srtd}
\end{table}

\subsection{Time-to-Solve Comparisons}\label{subsec4:4_time}
From a pure iteration count standpoint, SRTD and EVSS are largely incomparable. The EVSS formulation requires Newton's method to be run on a large system, while SRTD itself is iterative, and each SRTD iterate requires solving a smaller nonlinear system (meaning Newton's method must be used during Stage 1 of every SRTD iterate) and two linear systems. Hence, we include some information comparing the total solve time of the SRTD and EVSS methods for parameter values where both converge and are mesh stable. 

FEniCS is not native to Windows, so the following tests were performed on an Ubuntu VM running on a Windows machine via Windows Subsystem for Linux (WSL). The Ubuntu VM had access to 12 GB of RAM and 4 cores of an AMD Ryzen 5 3600 processor. The values listed in Tables \ref{tab:time_ldc_l1=0.01} and \ref{tab:time_jb_l1=0.01} below are the average time in seconds over three solves required for each method to converge. 

\begin{table}[htbp!]
\centering
\begin{tabular}{r|c|c|c}
\hline
Method\textbackslash $h$ & 2.50e-02 & 1.25e-02 & 6.25e-03 \\ \hline
SRTD             & 3.94        & 22.27       & 135.75      \\
EVSS             & 4.66        & 32.41       & 227.26     \\
\hline
\end{tabular}
\quad\quad
\begin{tabular}{r|c|c|c}
\hline
Method\textbackslash $h$ & 2.50e-02 & 1.25e-02 & 6.25e-03 \\
\hline
SRTD                                         & 3.90     & 22.31    & 135.52   \\
EVSS                                         & 4.73     & 32.50    & 228.10  \\
\hline
\end{tabular}
\caption{The average amount of time, in seconds, over three solves required by each method to solve the lid-driven cavity problem with $\lambda_{1}=0.01$ on the UCM model $\mu_{1}=\lambda_{1}$ (left) and the corotational Maxwell model $\mu_{1}=0$ (right) for various mesh sizes.}
\label{tab:time_ldc_l1=0.01}
\end{table}

\begin{table}[htbp!]
\centering
\begin{tabular}{r|c|c|c}
\hline
Method\textbackslash $h$ & 5.00e-02 & 2.50e-02 & 1.25e-02 \\ \hline
SRTD        & 2.61     & 12.86    & 73.90    \\
EVSS        & 2.78     & 17.92    & 124.70  \\
\hline
\end{tabular}
\quad\quad
\begin{tabular}{r|c|c|c}
\hline
Method\textbackslash $h$ & 5.00e-02 & 2.50e-02 & 1.25e-02 \\ \hline
SRTD        & 3.11     & 15.45    & 88.69    \\
EVSS        & 2.84     & 18.09    & 125.50  \\
\hline
\end{tabular}
\caption{The average amount of time, in seconds, over three solves required by each method to solve the journal-bearing problem with $\lambda_{1}=0.01$ on the UCM model $\mu_{1}=\lambda_{1}$ (left) and the corotational Maxwell model $\mu_{1}=0$ (right) for various mesh sizes.}
\label{tab:time_jb_l1=0.01}
\end{table}

Within the parameter range where both SRTD and EVSS converge and are mesh stable, SRTD reaches a solution within the desired relative tolerance of $10^{-9}$ noticeably faster than EVSS. The difference appears to grow as the mesh is refined.

\subsection{Newtonian Compared with Oldroyd 3}\label{subsec4:5_newt_vs_o3}
\begin{figure}[htbp!]
    \centering
    \begin{subfigure}{.5\textwidth}
        \centering
        \includegraphics[width=\linewidth]{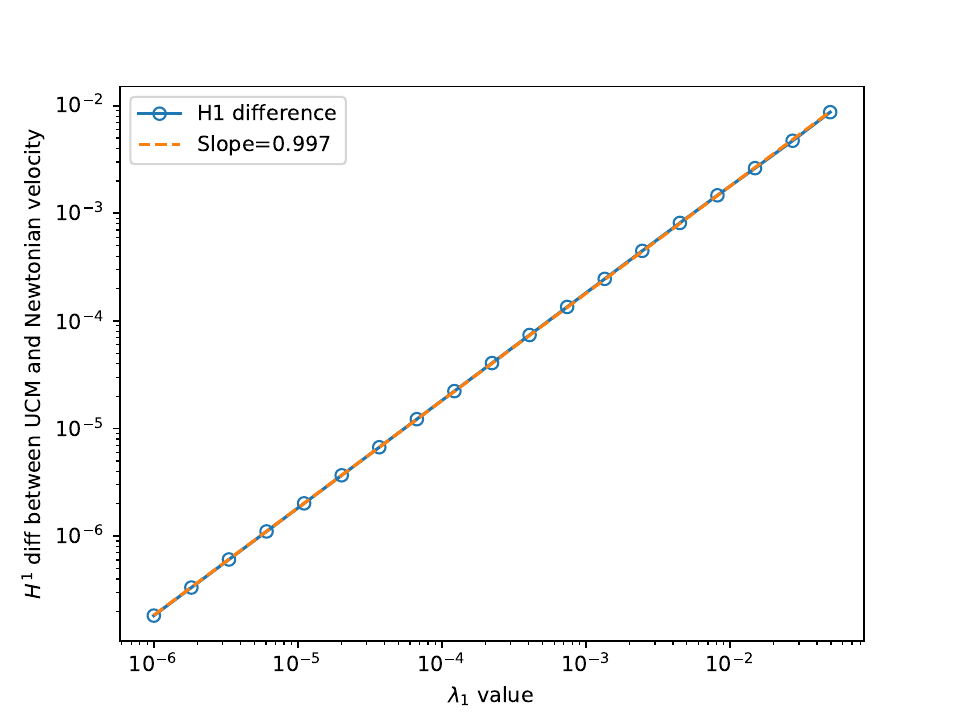}
    \end{subfigure}%
    \begin{subfigure}{.5\textwidth}
        \centering
        \includegraphics[width=\linewidth]{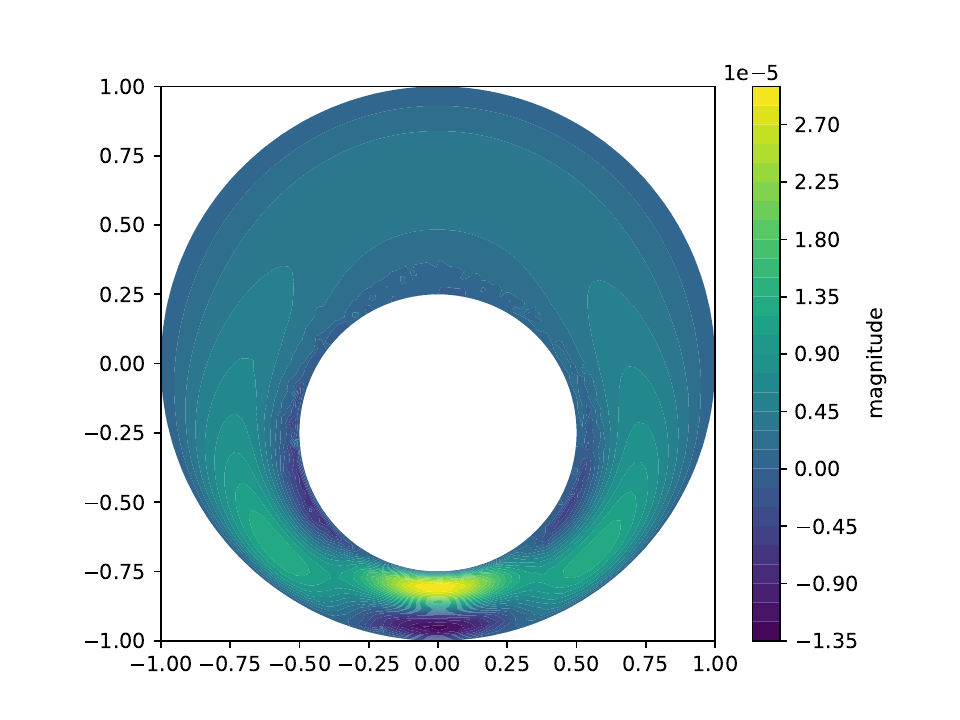}
    \end{subfigure}
    \caption{Left: the $H^{1}$ norm of the difference between the UCM solution, solved by the SRTD algorithm, and the Navier-Stokes/Newtonian solution for the journal-bearing problem, plotted against and increasing $\lambda_{1}$ values. Right: the difference between the azimuthal components (in the direction of increasing $\theta$ relative to the center of the bearing) of velocity for the UCM solution with $\lambda_{1}=0.01$ and Navier-Stokes/Newtonian solution.}
    \label{fig:compare_nse_jb}
\end{figure}

\begin{figure}[htbp!]
    \centering
    \begin{subfigure}{.5\textwidth}
        \centering
        \includegraphics[width=\linewidth]{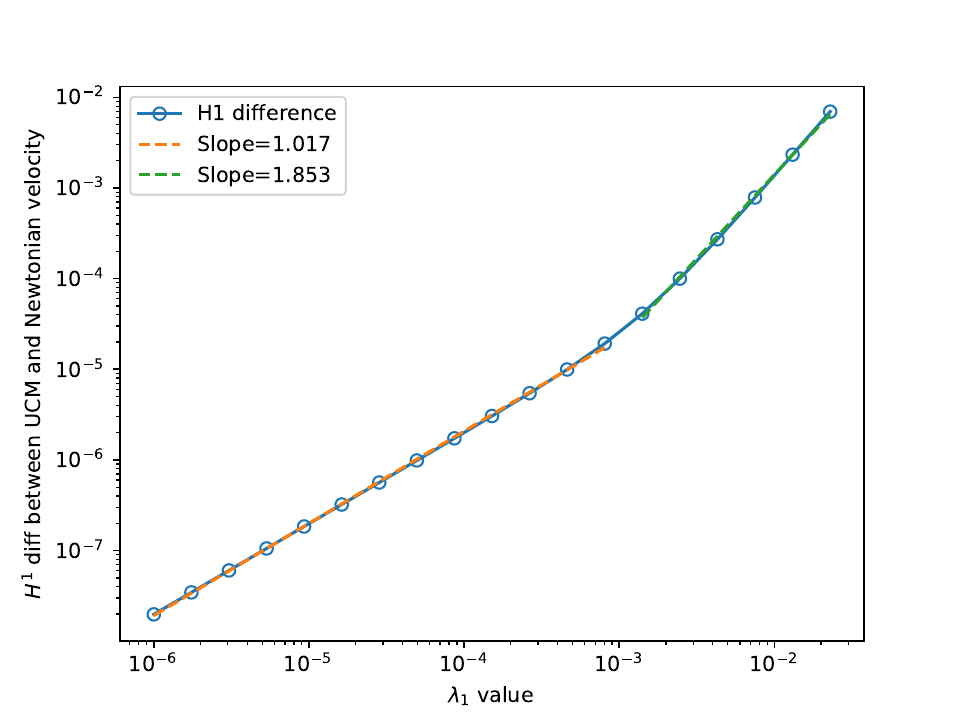}
    \end{subfigure}%
    \begin{subfigure}{.5\textwidth}
        \centering
        \includegraphics[width=\linewidth]{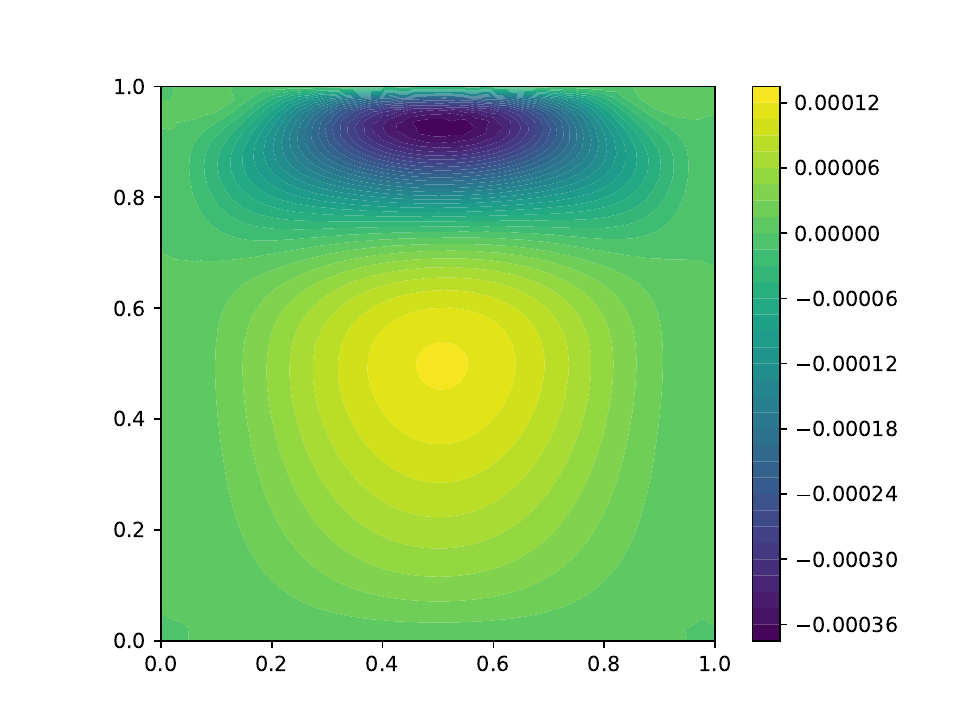}
    \end{subfigure}
    \caption{Left: the $H^{1}$ norm of the difference between the UCM solution, solved by the SRTD algorithm, and the Navier-Stokes/Newtonian solution for the lid-driven cavity problem, plotted against and increasing $\lambda_{1}$ values. Right: the difference between the horizontal components of velocity for the UCM solution with $\lambda_{1}=0.01$ and Navier-Stokes/Newtonian solution.}
    \label{fig:compare_nse_ldc}
\end{figure}

Given that the SRTD method (and EVSS to a lesser extent) failed to converge for even moderately-sized Weissenberg numbers, it is a valid question to ask whether the range of parameters in which SRTD does converge corresponds to a non-Newtonian solution which is substantially qualitatively different from the Newtonian solution. Observe from Figures \ref{fig:compare_nse_jb} (left) and \ref{fig:compare_nse_ldc} (left) the $H^{1}$ difference between the Navier-Stokes solution and the UCM solution. For nearly the entire of range of $\lambda_{1}$ values on the journal-bearing problem for which SRTD converges, the UCM solution diverges away from the Navier-Stokes solution linearly with $\lambda_{1}$. For the lid-driven cavity problem, the difference grows linearly with $\lambda_{1}$ until a certain point around $\lambda_{1}=1e-3$, when the relationship becomes nearly quadratic. 

\section{Conclusions}\label{sec5:conclusions}
The SRTD formulation and corresponding algorithm appears to be stable when the iterative process converges. It delivers comparable results relatively quickly compared to EVSS. By decoupling the problem, SRTD is capable of solving on more refined meshes than EVSS if the amount of computational resources is a limiting factor. SRTD was, in fact, still able to converge on more refined meshes than were included in this paper. This makes sense given the degrees of freedom data in Tables \ref{tab:jb_dof}, \ref{tab:ldc2d_dof}, and \ref{tab:ldc3d_dof}. Given that the number of degrees of freedom grows even more quickly with decreasing mesh size in three dimensions, this seems to make SRTD a suitable choice for solving three-dimensional flows. The drawback is that SRTD does not seem to be able to reach as high a Weissenberg number as EVSS.

\printbibliography

\end{document}